\documentclass[]{interact}
\usepackage{environ} 
\usepackage[dvipsnames]{xcolor}
\usepackage{hyperref}
\usepackage{graphicx}
\usepackage[normalem]{ulem}
\graphicspath{{./gfx/}}

\usepackage{epstopdf}

\usepackage[font=small,labelfont=bf]{caption}
\usepackage{subcaption}
\DeclareCaptionFormat{sqa}{#1~#3}
\usepackage{natbib}
\bibpunct[, ]{(}{)}{;}{a}{}{,}

\usepackage{tikz}
\usetikzlibrary{arrows.meta}
\usetikzlibrary{calc}
\usetikzlibrary{patterns}

\usepackage{mathrsfs}
\usepackage{bbm}
\usepackage{bm}

\newcommand{\Pb}{{\mathsf{P}}} 
\newcommand{\Qb}{{\mathsf{Q}}} 
\newcommand{\Eb}{{\mathsf{E}}} 
\newcommand{\ESS}{{\mathsf{ESS}}}
\newcommand{\Iin}{{\mathsf{I}_{\mrm{in}}}}
\DeclareMathOperator{\ARLFA}{\mathsf{ARLFA}}   
\newcommand{\ESADD}{{\mathsf{ESEDD}}}

\newcommand{\mrm}[1]{\mathrm{#1}}
\newcommand{\drm}{{\mrm{d}}}

\newcommand{\mc}[1]{\mathcal{#1}} 

\newcommand{\Nc}{{\mc{N}}}

\newcommand{\Lc}{{\mc{L}}}

\newcommand{\ccP}{\mathcal{P}}

\newcommand{\mbs}[1]{\bm{#1}} 
\newcommand{\alphab}{{\mbs{\alpha}}}
\newcommand{\betab}{{\mbs{\beta}}}
   
\newcommand{\pib}{{\mbs{\pi}}}

\newcommand{\mbb}[1]{\mathbb{#1}} 
\newcommand{\Rbb}{\mbb{R}} 
\newcommand{\Zbb}{\mbb{Z}} 

\newcommand{\class}{{\mbb{C}}}

\newcommand{\classPal}{{\mbb{C}_{\ccP}}(\alpha_0,\alpha_1)}

\newcommand{\Nbb}{\mathbb{N}}

\newcommand{\mb}[1]{\mathbf{#1}} 

\newcommand{\Xb}{{\mb{X}}}

\newcommand{\Bb}{\mb{B}}

\newcommand{\set}[1]{\left\{#1\right\}}

\newcommand{\brc}[1]{\left(#1\right)}
\newcommand{\brcs}[1]{\left[#1\right]}

\newcommand{\esssup}{\operatornamewithlimits{ess\,sup}}

\newcommand{\ignore}[1]{} 

\newcommand{\Fc}{{ \mathscr{F}}} 

\newcommand{\Bc}{{ \mathscr{B}}}

\newcommand{\cN}{{ \mathscr{N}}}

\newcommand{\alpham}{{\alpha_{\rm max}}}

\NewEnviron{notforprint*}{\BODY}
\NewEnviron{trivial}{\BODY}
\NewEnviron{ignore*}{}
\NewEnviron{skipproof}{}
\NewEnviron{noignore*}{\BODY}

\newcommand{\nolabel}[1]{}

\hypersetup{
    colorlinks=true,
    bookmarksnumbered=true,
    bookmarksopen=true,
    linkcolor=blue,
    citecolor=blue,
    urlcolor=blue,
    unicode=true,
    breaklinks=true
}
\renewcommand{\geq}{\geqslant}
\renewcommand{\leq}{\leqslant}
\renewcommand{\ge}{\geqslant}
\renewcommand{\le}{\leqslant}

\renewcommand{\cup}{\,{\textstyle\bigcup}\,}

\tikzset{>={Latex[length=8pt, width=4pt]}}

\DeclareMathOperator{\Hyp}{\mathcal{H}}

\newcommand{\cB}{\mathcal{B}}

\newcommand{\cE}{\mathcal{E}}

\newcommand{\cL}{\mathcal{L}}

\newcommand{\cP}{\mathcal{P}}

\theoremstyle{plain}

\newtheorem*{lemma*}{Lemma}
\newtheorem{theorem}{Theorem}
\newtheorem*{theorem*}{Theorem}

\newtheorem*{corollary*}{Corollary}

\newtheorem*{proposition*}{Proposition}
\theoremstyle{remark}

\newtheorem*{assumption*}{Assumption}
\theoremstyle{definition}

\newtheorem*{definition*}{Definition}

\DeclareFontFamily{U}{matha}{\hyphenchar\font45}
\DeclareFontShape{U}{matha}{m}{n}{
      <5> <6> <7> <8> <9> <10> gen * matha
      <10.95> matha10 <12> <14.4> <17.28> <20.74> <24.88> matha12
      }{}
\DeclareSymbolFont{matha}{U}{matha}{m}{n}
\DeclareMathSymbol{\abscont}{3}{matha}{"21}

\newcommand{\um}{\underline{K}}
\newcommand{\om}{\overline{K}}

\newcommand{\whdelta}{{\widehat{\delta}}}
\newcommand{\whT}{{\widehat{T}}}

\newcommand{\whd}{{\widehat{d}}}

\newcommand{\dotb}{\overset{\bm .}}
\newcommand{\var}{{\mathsf{var}}} 

\NewEnviron{revone*}{{\color{Magenta}\BODY~}}
\NewEnviron{revtwo*}{{\color{OliveGreen}\BODY~}}
\NewEnviron{revmisc*}{{\color{ProcessBlue}\BODY~}}

\begin{document}

\title{Beyond boundaries: Gary Lorden's groundbreaking contributions to sequential analysis}

\author{
    \name{Jay Bartroff\textsuperscript{a} and Alexander~G. Tartakovsky\textsuperscript{b} 
    \thanks{CONTACT: Jay Bartroff. Email: bartroff@austin.utexas.edu} 
}
    \affil{%
        \textsuperscript{a} University of Texas at Austin, Austin, Texas, USA; 
        \textsuperscript{b}AGT StatConsult, Los Angeles, California, USA}
}

\maketitle

\begin{abstract}
Gary Lorden provided several fundamental and novel insights into sequential hypothesis testing and changepoint detection. In this article, we provide an overview of 
Lorden's contributions in the context of existing results in those areas, and some extensions made possible by Lorden's work.  We also mention some of Lorden's significant consulting work, including as an expert witness and for NASA, the entertainment industry, and Major League Baseball.
\end{abstract}


\section{Introduction} \label{sec:intro}

The purpose of this article is to provide an overview of Gary Lorden's  significant  contributions to the field of sequential analysis. But first, to give a sense of Lorden's remarkable life and personality, in Sections~\ref{sec:bio}-\ref{sec:consulting} we give a brief biography of Lorden and highlight some of his extra-academic work.

\begin{figure}
    \centering
    \includegraphics[scale = .14]{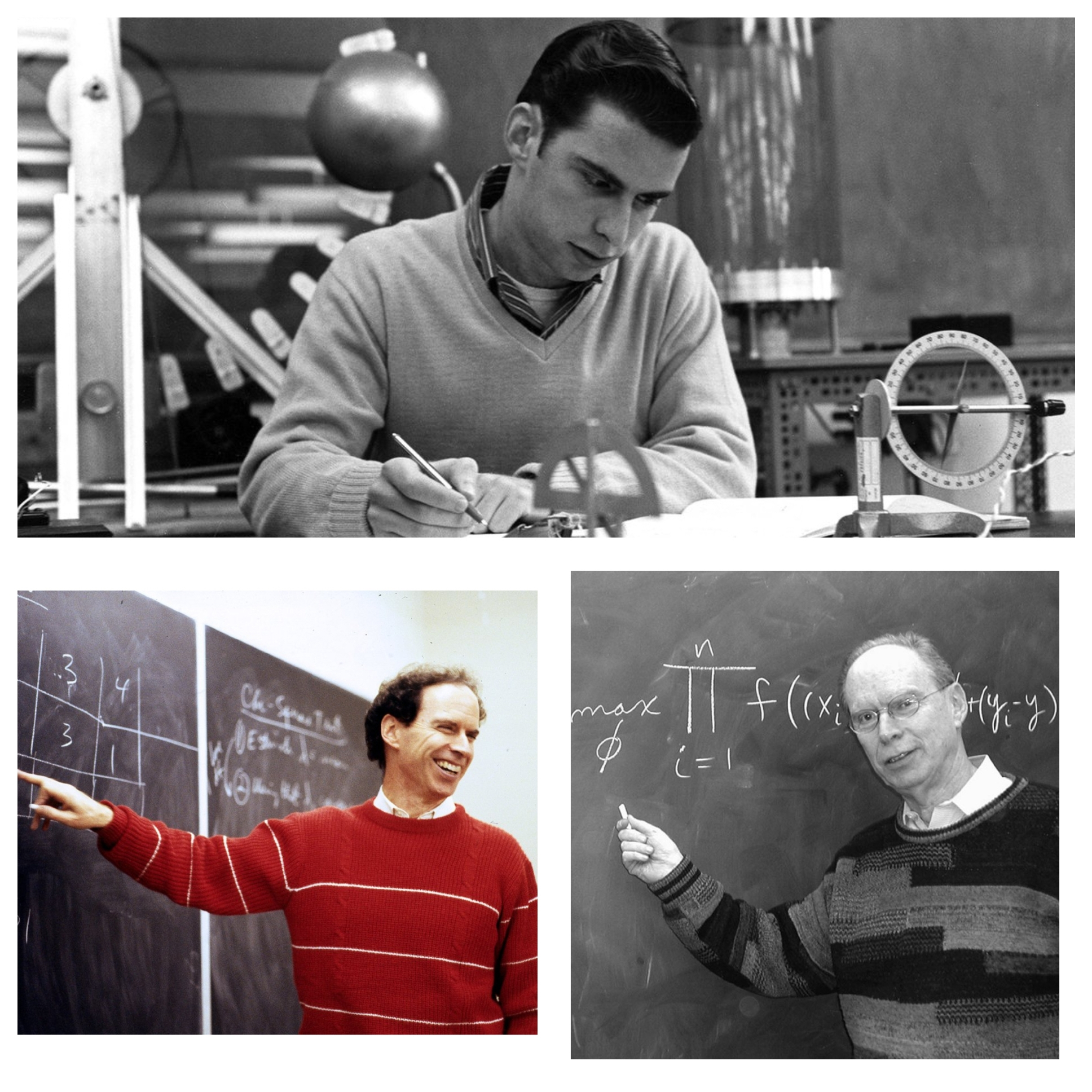}
    \caption{Counter-clockwise from top: Gary Lorden as an undergraduate at Caltech in 1959; lecturing in 1987; and lecturing circa 2010.  Photos courtesy of Caltech.}
    \label{fig:lorden_pics}
\end{figure}

\subsection{Biographical Sketch}\label{sec:bio}

Gary Allen Lorden was born in Los Angeles, California, on June 10, 1941. Lorden entered Caltech as a freshman in 1958, and Lorden's undergraduate contemporaries at Caltech included a number of others who would also go on to be notable statisticians including Larry Brown, Peter Bickel, Brad Efron, and Carl Morris. Lorden received a BS from Caltech in 1962 and a PhD from Cornell University in 1966 under the supervision of Jack Kiefer. After a faculty position at Northwestern University, Lorden rejoined Caltech in 1968, where he stayed until his retirement in 2009.

Beyond research and teaching, Lorden was known for his leadership roles at Caltech. He served as dean of students from 1984 to 1988, vice president for student affairs from 1989 to 1998, and acting vice president for student affairs in 2002. He was executive officer (Caltech's version of department chair) for the mathematics department from 2003 to 2006.

Gary and his wife, Louise, were both accomplished pianists and enjoyed playing and singing duets for guests, often students, at their home in Pasadena.  They also enjoyed showing students how well they danced together. Lorden also liked to act and regularly participated in plays put on at Caltech.

Lorden passed away on October 25, 2023 at the age of 82.  He is survived by his children Lisa and Diana, and his wife Louise passed away in 2015.

\subsection{Consulting, Hollywood, and Major League Baseball}\label{sec:consulting}

Lorden was known for his creativity and generosity of ideas, so  it is not surprising that throughout his career he was in demand as an academic collaborator and consultant, and much of his early work on changepoint detection arose from collaborations with researchers at the Jet Propulsion Laboratory (JPL) in Pasadena, California, which is managed for NASA by Lorden's home institution of Caltech.  

But Lorden was also an uncommonly effective, engaging, and entertaining communicator and this caused many outside academics to seek him out as a consultant too. Lorden would routinely be interviewed by reporters and appear on the evening news, explaining things like the odds of winning the latest lottery jackpot. Lorden also ``moonlighted'' as an expert witness in court cases involving statistics and math.  

Eventually Hollywood came calling, and in 2005 Lorden was asked to be the math consultant for a (then) new  CBS television show called \textit{NUMB3RS}.  In a bit of art imitating life, the show was about a Caltech professor who used math to help the FBI solve crimes, and Lorden worked with the show's writers to accurately incorporate mathematical topics into the storylines. An aspect of the relationship that Lorden particulalry relished, that he related to us, was that the show's creators initially considered setting the show at MIT (Caltech's rival), but decided to change the venue after learning more about Caltech and Lorden. The show would go on to be a hit, running for 118 episodes over 6 seasons.  If a viewer of the show knew something of Lorden's work they could see some of his favorite topics (and many topics discussed below in this article) woven into the episodes' plot lines including changepoint detection, hypothesis testing, Bayesian methods, gambling math, cryptography, and sports statistics among others. With Keith Devlin, Lorden wrote a popular general audience book \citep{Devlin07} on the mathematical topics appearing in the show.  For example,  \citet[][Chapter~4]{Devlin07} explain the basic concept of 
changepoint detection: 
\begin{quote}
    ``the determination that a definite change has occurred, as opposed to normal fluctuations,''
\end{quote} and goes on to discuss this method's importance for quick response to potential bioterrorist attacks and for designing efficient algorithms to pinpoint various kinds of criminal activity, such as to detect an increase in crime rates in certain geographical areas and to track changes in financial transactions that could be criminal.

Another high-profile consulting project came in 2018 when Lorden, Bartroff, and others were chosen by the Commissioner of Major League Baseball (MLB)  to be statisticians on a committee with physicists, engineers, and baseball experts studying MLB’s then-recent surge in home runs.  The committee studied vast amounts of Statcast game data,  as well as laboratory tests on the properties of the baseball that can affect home run production, including Lorden traveling to Costa Rica to inspect baseball manufacturer Rawlings' production plant there. The committee made recommendations~\citep{Nathan18} to MLB and Rawlings for future monitoring, testing, and storage of baseballs.

\subsection{The Remainder of this Article}

The remainder of this article is devoted to describing Lorden's fundamental and far-reaching results in 
sequential hypothesis testing and changepoint detection, which we aim to present within the context of those areas. Beginning with hypothesis testing in Section~\ref{sec:SHT}, after describing the testing setup and optimality of the sequential probability ratio test (SPRT) in Section~\ref{ssec:PF}, Lorden's findings on multi-parameter testing and their application to near-optimality of the multihypothesis SPRT are covered in Section~\ref{ssec:MSPRT}. We then cover Lorden's fundamental inequality for excess over the boundary in Section~\ref{sec:lord.ineq}.  Additionally, we explore Lorden's contributions to the Keifer-Weiss problem of testing while minimizing the expected sample size at a parameter value between the hypotheses and other results stemming from his work (Section~\ref{ssec:2SPRT}), optimal testing of composite hypotheses (Section~\ref{ssec:GLR}), and optimal multistage testing (Section~\ref{sec:multistage}). In Section~\ref{sec:CPD} we cover Lorden's fundamental minimax changepoint detection theory and related advancements in the field.

\section{Sequential Hypothesis Testing}\label{sec:SHT}

In this section, we delve into Lorden's contributions to hypothesis testing, encompassing the challenges of excess over the boundaries of random walks, the formulation of nearly optimal multi-decision sequential rules (including multihypothesis tests and multistage tests), and the modified Keifer--Weiss problem. These topics were extensively explored in the seminal Lorden's papers \citep{Lorden1,lorden-ams70,lorden-as76,Lorden-AS77,Lorden-ams72,Lorden-AS73,lorden-ptrf80,Lorden83}.

\subsection{Multihypothesis Testing}
\subsubsection{The General Multihypothesis Testing Problem} \label{ssec:PF}

One of Lorden's major fundamental contributions is the proposal of a multihypothesis sequential test that achieves third-order asymptotic optimality in the i.i.d.\ case, 
akin to Wald's SPRT when error probabilities are small. Further details will be provided in Subsection~\ref{ssec:MSPRT}. Additionally, extensions to a more general non-i.i.d. case, where observations may exhibit dependency and non-identical distributions, have been explored in \citet{Lai-as81-SPRT,TartakovskySISP98,Tartakovsky_book2020,TartakovskyAMSA2024,TNB_book2014}.

We begin with formulating the following multihypothesis testing problem as addressed by \citet{Lorden-AS77}. Let $(\Omega, \Fc ,\Fc_n, \Pb )$, 
$n \in \Zbb_+ = \{0,1, 2,\ldots \}$, be a filtered probability space, where the sub-$\sigma$-algebra $\Fc_n=\sigma(\Xb^n)$ of~$\Fc$ is generated by the sequence of random variables 
$\Xb^n = \{X_t, \: 1 \le t \le n\}$ observed up to time~$n$, which is defined on the space~$(\Omega, \Fc)$ ($\Fc_0$ is trivial). 
The focus lies on the $N$-decision problem of testing the hypotheses $\Hyp_i:~ \Pb = \Pb_i$, $i = 1,\dots,N$, where $\Pb_1,\dots,\Pb_N$ are given probability measures 
assumed to be locally mutually 
absolutely continuous, i.e., their restrictions $\Pb_i^{n}=\Pb|_{\Fc_n}$ and~$\Pb_j^{n}= \Pb_j |_{\Fc_n}$ to $\Fc_n$ are equivalent for all $1 \le n< \infty$ and 
all $i, j =1,\dots,N$, $i \neq j$. Let
$\Qb^{n}$ be a restriction to $\Fc_n$ of a non-degenerate $\sigma$-finite measure $Q$ on $(\Omega, \Fc)$. 

Assume that the observed random variables $X_1,X_2,\dots$ are independent and identically distributed (i.i.d.), so that under $\Pb_i$ the sample $\Xb^n = (X_1,\dots,X_n)$ has joint density~$p_{i}(\Xb^n)$ with respect to the dominating measure $\Qb^{n}$ for all $n\in \Nbb$, which can be expressed as
\begin{equation}\label{Jointdensiid}
p_{i}(\Xb^n) = \prod_{t=1}^n f_{i}(X_t), \quad i=1,\dots,N,
\end{equation}
where $ f_{i}(X_t)$ represents the respective density for $X_t$ under the hypothesis $\Hyp_i$. Hence, the interest lies in determining which of $N$ given 
densities $f_1, \dots,f_N$ is true.

A multihypothesis sequential test is a pair~$\delta=(T,d)$, 
where $T$ is a stopping time with respect to the filtration $\{\Fc_n \}_{n \in \Zbb_+}$ and $d=d(\Xb^{T})$ is an $\Fc_T$-measurable terminal decision function 
with values in the set $\{1,\dots,N\}$. Specifically, $d=i$ means that the hypothesis~$\Hyp_i$ is accepted upon stopping,  
$\set{d=i}=\set{T<\infty, ~ \delta ~ \text{accepts~$\Hyp_i$}}$. Let $\alpha_{ij}(\delta)=\Pb_i(d=j)$, $i \neq j$, $i,j=1,\dots,N$, denote the error probabilities of 
the test~$\delta$, i.e., the probabilities of accepting the hypothesis~$\Hyp_j$ when $\Hyp_i$ is true. 

Introduce the class of tests with probabilities of errors $\alpha_{ij}(\delta) $ that do not exceed the prespecified numbers $0<\alpha_{ij}<1$:
\begin{equation}\label{Mclasses}
\class(\alphab)  = \set{\delta: \alpha_{ij}(\delta)  \le \alpha_{ij} ~\text{for} ~ i, j = 1,\dots,N, \, i \neq j},
\end{equation}
where $\alphab=(\alpha_{ij})$ is a matrix of given error probabilities that are positive numbers less than~$1$ (the diagonal entries $\alpha_{ii}$ are immaterial).

Let $\Eb_i$ denote the expectation under the hypothesis~$\Hyp_i$ (i.e., under the measure $\Pb_i$). The objective is to discover a sequential test that 
would minimize the expected sample sizes $\Eb_i [T]$ for all hypotheses $\Hyp_i$, $i=1,\dots,N$, at least approximately.

For $n\in \Nbb$, define the likelihood ratio (LR) and the log-likelihood ratio (LLR) processes between the hypotheses $\Hyp_i$ and~$\Hyp_j$ as
\begin{align*}
\Lambda_{ij}(n) & = \frac{\drm \Pb_i^{n}}{\drm \Pb_j^{n}}(\Xb^n) = \frac{p_{i}(\Xb^n)}{p_{j}(\Xb^n)} = \prod_{t=1}^n \frac{f_{i}(X_t)}{f_{j}(X_t)},
\\
\lambda_{ij}(n) & = \log \Lambda_{ij}(n) = \sum_{t=1}^n \log \brcs{\frac{f_{i}(X_t)}{f_{j}(X_t)}} .
\end{align*}


In a particular case of two hypotheses $\Hyp_1$ and $\Hyp_2$ ($N=2$), \citet{wald45,wald47} introduced 
the {\em Sequential Probability Ratio Test} (SPRT). 
Let $Z_t = \log[f_1(X_t)/f_2(X_t)]$ be the LLR for the observation~$X_t$, so the LLR for the sample $\Xb^n$ 
is the sum 
\begin{equation}\label{llr.sprt}
\lambda_{12}(n) = \lambda_n=  \sum_{t=1}^n Z_t, \quad n = 1, 2,\dots
\end{equation}
Letting $a_0<0$ and $a_1 > 0$ be thresholds, Wald's SPRT $\delta_*(a_0,a_1) = (T_*,d_*)$ is 
\begin{equation}\label{SPRT}
 T_*(a_0,a_1) = \inf\set{n \ge 1: \lambda_n \notin (a_0, a_1)}, \quad
d_*(a_0,a_1) = \begin{cases}
1 & \text{if}~~ \lambda_{T_*} \ge a_1\\
2 & \text{if}~~ \lambda_{T_*} \le a_0 .
\end{cases}
\end{equation}

In the case of two hypotheses, the class of tests \eqref{Mclasses} is defined as
$$
\class(\alpha_0, \alpha_1)  = \set{\delta: \alpha_{0}(\delta)  \le \alpha_{0} ~\text{and } ~ \alpha_{1}(\delta)  \le \alpha_{1}},
$$
where $\alpha_0$ represents the upper bound on the Type I error (false positive) probability $\alpha_0(\delta) = \alpha_{12}(\delta)$, and $\alpha_1$ represents
the upper bound on the Type II error (false negative)  probability $\alpha_1(\delta) = \alpha_{21}(\delta)$.  

Wald's SPRT possesses  an extraordinary optimality property: it minimizes both the expected sample sizes $\Eb_1 [T]$ and~$\Eb_2 [T]$ within the class of sequential 
(and non-sequential) tests $\class(\alpha_0, \alpha_1)$  as long as the observations are i.i.d.\ under both hypotheses. In other words:
$\Eb_i[T_*] = \inf_{\delta\in \class(\alpha_0,\alpha_1)} \Eb_i[T]$ for $i=0, 1$,
as established by  \citet{wald48} through a Bayesian approach.

\citet{Lai-as81-SPRT} proved that the SPRT is also first-order asymptotically optimal as $ \max(\alpha_0,\alpha_1) \to 0$ for general non-i.i.d.\ models 
with dependent and non-identically 
distributed observations when the normalized log-likelihood ratio $n^{-1} \lambda_n$ converges $1$-quickly to finite numbers $I_i$ under $\Pb_i$.

The central idea of Lorden's investigation, elaborated in detail in Section~\ref{ssec:MSPRT}, is that, similar to how the SPRT is strictly optimal in the class $\class(\alpha_0,\alpha_1)$ for any error probabilities $\alpha_0$ and $\alpha_1$, combinations of SPRTs exhibit third-order asymptotic optimality for multihypothesis testing problems involving any finite number of densities when probabilities of errors are small.

\subsubsection{Near Optimality of the Multihypothesis SPRT}\label{ssec:MSPRT}

The problem of sequentially testing many hypotheses is 
substantially more complex than that of testing two hypotheses. Identifying an optimal test in the class~\eqref{Mclasses} that minimizes expected sample sizes for all hypotheses $\Hyp_1, \dots, \Hyp_N$, is daunting.  Hence, a significant portion of the development of sequential multihypothesis testing in the 20th century has focused on 
the exploration of certain combinations of one-sided SPRTs. See \citet{Armitage,Cher,KieferSacks-AMS1963,Lorden1,Lorden-AS77}. 

The results of Lorden's ingenious paper \citet{Lorden-AS77} are of fundamental importance as they establish third-order asymptotic optimality of the 
accepting multihypothesis test that
he proposed. More specifically, Lorden established that just as the SPRT is optimal in the class $\class(\alpha_0,\alpha_1)$ for testing two hypotheses, 
certain combinations of one-sided SPRTs are nearly optimal in a third-order sense in the class $\class(\alphab)$, i.e., subject to error probability constraints expected sample sizes are
minimized to within the negligible additive $o(1)$ term:
\begin{equation}\label{TOasymptOpt}
\inf_{\delta\in\class(\alphab)} \Eb_i[T] = \Eb_i[T_*] + o(1) \quad \text{as}~ \alpham\to 0 \quad \text{for all} ~ i=1,\dots,N,
\end{equation}
where $\alpham = \max_{1\le i,j \le N, i \neq j} \alpha_{ij}$ and $T_*$ is the stopping time of the multihypothesis test $\delta_*$, which is defined below.

We now define a test proposed by Lorden, which we will refer to as the accepting Matrix SPRT.
Write $\cN=\{1,\dots,N\}$. For a threshold matrix $\mathbf{A}=(A_{ij})_{i,j\in\cN}$, with $A_{ij} > 0$ for $i\ne j$ and the $A_{ii}$ are immaterial ($0$, say), define 
the Matrix SPRT (MSPRT) $\delta_* = (T_*, d_*)$, built on one-sided SPRTs between the hypotheses 
$\Hyp_i$ and~$\Hyp_j$, as follows:
\begin{equation}\label{MSPRT1}
\text{Stop at the first $n \geq 1$ such that, for some $i$}, ~  \Lambda_{ij}(n) \geq A_{ji} ~  \text{for all $j \neq i$},
\end{equation}
and accept the unique~$\Hyp_i$ that satisfies these inequalities. Note that for $N=2$ the MSPRT coincides with Wald's SPRT.  

Let $a_{ji}= \log A_{ji}$. Introducing the Markov accepting times for the hypotheses~$\Hyp_i$ as
\begin{equation} \label{taui}
T_i=\inf\set{ n \geq 1: \min_{\stackrel{1 \le j \le N}{ j \neq i}} \brcs{\lambda_{ij}(n) -a_{ji}} \geq 0 } , \quad i =1,\dots,N,
\end{equation}
the test in \eqref{MSPRT1} can also be written in the following form:
\begin{equation} \label{D1}
T_*=\min_{1\le j \le N} T_j, \qquad d_*=i \quad \mbox{if} \quad T_*= T_i.
\end{equation}
Thus, in the MSPRT, each component SPRT is extended until, for some $i\in \cN$, all $N-1$ SPRTs involving $\Hyp_i$ accept~$\Hyp_i$.
 
The MSPRT is not strictly optimal for $N>2$ but it is a good approximation to the  optimal multihypothesis test. Under certain conditions and with some 
choice of the threshold matrix~$\mathbf{A}$, it minimizes the expected sample sizes~$\Eb_i [T]$ for all $i=1,\dots, N$ to within a vanishing $o(1)$ term  
for small error probabilities; see \eqref{TOasymptOpt}.

Consider first the first-order asymptotic criterion:
Find a multihypothesis test $\delta_*(\alphab)=(d_*(\alphab),T_*(\alphab))$ such that 
\begin{equation} \label{FOASoptiM}
\lim_{\alpham\to0} \frac{\inf_{\delta\in \class(\alphab)}\Eb_i[T]}{\Eb_i[T_*(\alphab)]} =1 \quad \text{for all}~i=1,\dots,N.
\end{equation}

Using Wald's likelihood ratio identity, it is easily shown that $\alpha_{ij}(\delta_*) \le \exp(-a_{ij})$ for $i,j = 1,\dots,N$, $i \neq  j$, so selecting
$a_{ji}  = |\log \alpha_{ji}|$ implies $\delta_* \in \class(\alphab)$. These inequalities are similar to Wald's in the binary hypothesis case and are very imprecise.
Using Wald's approach it is rather easy to prove that the MSPRT with boundaries $a_{ji}  = |\log \alpha_{ji}|$ is first-order asymptotically optimal, minimizing expected sample sizes
as long as the Kullback-Leibler information numbers  $I_{ij}=\Eb_i [\lambda_{ij}(1)]$  are positive and finite; see \citet[][Section 4.3.1]{TNB_book2014}. 

In his ingenious paper, \citet{Lorden-AS77} substantially improved this result showing that with a sophisticated design that includes accurate estimation of thresholds 
accounting for overshoots, the MSPRT is nearly optimal in the third-order sense~\eqref{TOasymptOpt}. 

Specifically, assume the second-moment condition
\begin{equation}\label{2ndmomentcond}
\Eb_i [\lambda_{ij}(1)]^2 < \infty, \quad i,j = 1,\dots,N 
\end{equation}
and define the numbers
\begin{equation}\label{Upsilonij}
\cL_{ij}= \exp\set{- \sum_{n=1}^\infty \frac{1}{n}  [ \Pb_j (\lambda_{ij}(n) >0)  + \Pb_i(\lambda_{ij}(n) \le 0)]}, \quad i,j=1,\dots,N.
\end{equation}
These numbers are symmetric, $\cL_{ij}=\cL_{ji}$, and $0<\cL_{ij} \le 1$ ($\cL_{ii}\equiv 1$). 
Furthermore,  $\cL_{ij} = 1$ only if the measures $\Pb_i^n$ and~$\Pb_j^n$ are singular so that the absolute continuity assumption is violated. 

For $i,j\in \cN$ ($i\neq j$) and $a>0$, define one-sided SPRTs
\begin{equation}\label{onesidedij}
\tau_{ij}(a) = \inf\set{n \geq 0: \lambda_{ij}(n) \geq a} .
\end{equation}
Using a renewal-theoretic argument,  the numbers $\cL_{ij}$ are tightly related to the overshoots in the one-sided tests. If the LLR $\lambda_{ij}(1)$ is 
non-arithmetic under $\Hyp_i$, then
\begin{equation} \label{zetaell}
 \cL_{ij}= \zeta_{ij} I_{ij}, \quad \zeta_{ij} = \lim_{a\to\infty} \Eb_i \set{\exp\brcs{-(\lambda_{ij}(\tau_{ij}(a))-a)}}
\end{equation}
(see, e.g., Theorem~3.1.3 in \citet{TNB_book2014}).

It turns out that the $\cL$-numbers play a significant role both in the Bayes and the frequentist frameworks. They facilitate the adjustment of boundaries necessary 
to achieve optimality.

Consider the Bayes multihypothesis problem with the prior distribution of hypotheses $\pib=(\pi_0(1), \dots,\pi_0(N))$, where  $\pi_0(i)=\Pb(\Hyp_i)$, and the 
loss  incurred when stopping at time $T=n$ 
and making the decision $d=j$ while the hypothesis $\Hyp_i$ is true is $L_n(\Hyp_i, d=j, \Xb^n) = L_{ij} + c n$, where $c>0$ is the cost of making one observation or 
sampling cost and where $0<L_{ij} < \infty$ for $i \ne j$ and 0 if $i=j$.

 The average (integrated) risk of the test~$\delta=(T, d)$ is
\[
\rho_c^\pi(\delta) = \sum_{i=1}^N \pi_0(i) \brcs{\sum_{j=1}^N L_{ij} \Pb_i(d=j) + c \, \Eb_i [T]}.
\]
It follows from Theorem~1 of \citet{Lorden-AS77} that, as $c\to 0$, the MSPRT~$\delta_*$ defined in~\eqref{MSPRT1} with the thresholds $A_{ji}(c) = (\pi_0(j)/\pi_0(i)) L_{ji}\cL_{ij} / c$ 
is asymptotically third-order optimal (i.e., to within~$o(c)$) under the second moment condition~\eqref{2ndmomentcond}:
\[
\rho_c^\pi(\delta^*)= \inf_{\delta} ~\rho_c^\pi(\delta)  +o(c) \quad \text{as}~~c\to 0,
\]
where infimum is taken over all sequential or non-sequential tests.

Using this Bayes asymptotic optimality result, it can be proven that the MSPRT is also nearly optimal to within~$o(1)$ with respect to the expected sample sizes~$\Eb_i [T]$ for all hypotheses among all tests with constrained error probabilities. In other words, the MSPRT has an asymptotic property similar to the exact optimality of the SPRT for two hypotheses. This result is more practical than the above Bayes optimality.

The following theorem provides detailed specifications, resembling Theorem 4 and its corollary in \citet{Lorden-AS77}. Recall that $\alpha_{ij}(\delta) = \Pb_i(d=j)$ represents 
the probability to 
erroneously accept the hypothesis~$\Hyp_j$ when $\Hyp_i$ is true. In addition, denote as $\tilde{\alpha}_i(\delta)= \Pb_i(d\neq i)$ the probability of erroneously rejecting $\Hyp_i$ 
when it is true, and $\beta_j(\delta) =\sum_{i=1}^N w_{ij} \Pb_i(d=j)$ as the weighted probability of accepting~$\Hyp_j$, where $(w_{ij})_{i,j\in \Nc}$ is a given matrix of positive weights. 
Recall the definition of the class of tests~\eqref{Mclasses} for which the probabilities of errors $\Pb_i(d=j)$ do not exceed prescribed values $\alpha_{ij}$ and introduce two more classes 
that upper-bound the weighted probabilities of errors $\beta_j(\delta)$ and probabilities of errors $\tilde{\alpha}_i(\delta)$, respectively,
\begin{align}
\overline{\class}(\betab) & = \set{\delta: \beta_j(\delta) \leq \beta_j ~~ \text{for}~ j =1,\dots,N}, \label{Classbeta}
\\
\tilde{\class}(\tilde{\alphab}) & = \set{\delta: \tilde{\alpha}_i(\delta) \leq \tilde{\alpha}_i ~~ \text{for}~ i=1,\dots,N}. \label{Classalphatilde}
\end{align}

If $A_{ij}=A_{ij}(c)$ is a function of the small parameter~$c$, then the error probabilities $\alpha_{ij}^*(c)$, $\tilde{\alpha}^*_i(c)$ 
and~$\beta_j^*(c)$ of the MSPRT $\delta_*(c)$ are also functions of this parameter, and if $A_{ji}(c)\to\infty$, then $\alpha_{ij}^*(c), \beta^*_j(c) \to 0$ as $c\to0$. Note that
 $\tilde{\alpha}^*_i(c)=\sum_{j\neq i}\alpha^*_{ij}(c)$, so it also goes to zero as $c \to 0$.
We denote as $\betab^*(c)$ the vector~$(\beta_1^*(c), \dots,\beta_N^*(c))$, as $\tilde{\alphab}^*(c)$ the vector 
$(\tilde{\alpha}_1^*(c), \dots,\tilde{\alpha}_N^*(c))$ and as $\alphab^*(c)$ the matrix $(\alpha_{ij}^*(c))_{i,j\in \cN}$.

\begin{theorem}[MSPRT near optimality]\label{Th:MSPRTnearopt}
 Assume that the second moment condition~\eqref{2ndmomentcond} holds.
\begin{description}
\item [(i)] If the thresholds in the MSPRT are selected as $A_{ji}(c) = w_{ji} \cL_{ij}/c$, $i,j =1,\dots,N$, then
\begin{equation}\label{MSPRTAOiidbetaj}
\Eb_i [T^*(c)] = \inf_{\delta \in \overline{\class}(\betab^*(c))} \Eb_i [T]  + o(1) \quad \text{as}~ c \to 0 ~~ \text{for all}~~ i=1,\dots,N,
\end{equation}
i.e., the MSPRT minimizes to within~$o(1)$ the expected sample sizes among all tests whose weighted error probabilities are less than or equal to those of~$\delta^*(c)$.

\item [(ii)] For any matrix $\Bb=(B_{ij})$ ($B_{ij} >0$, $i \neq j$), let $A_{ji}=B_{ji}/c$. The MSPRT~$\delta^*(c)$ asymptotically minimizes the expected sample sizes for all 
hypotheses to within~$o(1)$ as $c\to0$ among all tests whose error probabilities~$\alpha_{ij}(\delta)$ are less than or equal to those of~$\delta^*(c)$ as well as whose error 
probabilities~$\tilde{\alpha}_{i}(\delta)$ are less than or equal to those of~$\delta^*(c)$, i.e.,
\begin{equation}\label{MSPRTAOiidalphaij}
\Eb_i [T^*(c)] = \inf_{\delta \in \class(\alphab^*(c))} \Eb_i [T]  + o(1) \quad \text{as}~~ c \to 0 \quad \text{for all}~~ i=1,\dots, N
\end{equation}
and
\begin{equation}\label{MSPRTAOiidalphai}
\Eb_i [T^*(c)] = \inf_{\delta \in \tilde{\class}(\tilde{\alphab}^*(c))} \Eb_i [T]  + o(1) \quad \text{as}~~ c \to 0 \quad \text{for all}~~ i=1,\dots, N.
\end{equation}
\end{description}
\end{theorem}

The intuition behind these results is that since the MSPRT is a combination of one-sided SPRTs $\tau_{ij}(a_{ji})$ defined in \eqref{onesidedij} and since the 
$\zeta_{ij}=\cL_{ij}/I_{ij}$ are correction factors to the error probability bound $\Pb_j(\tau_{ij}(a_{ji})<\infty) \le e^{- a_{ji}}$, the asymptotic approximation 
\[
\Pb_j(\tau_{ij}(a_{ji})<\infty) =  \zeta_{ij} e^{-a_{ji}} (1+o(1)) \quad \text{as} ~ a_{ji} \to \infty,
\]
works well even for moderate values of~$a_{ji}$. So taking $a_{ji}= \log (I_{ij}/\cL_{ij}\alpha)$ allows one to attain a nearly optimal solution in the frequentist problem. 
The proofs of these results are extremely tedious and require many non-standard and sophisticated mathematical tools developed by Lorden.

Notice that Theorem~\ref{Th:MSPRTnearopt} only addresses the asymptotically symmetric case where
\begin{equation}\label{symcasePrEr}
\lim_{c\to0} ~\frac{\log \beta_{j}^*(c)}{\log\beta_{k}^*(c)} =1, \quad \lim_{c\to0} ~\frac{\log \tilde{\alpha}_{i}^*(c)}{\log\tilde{\alpha}_{k}^*(c)} =1  \quad \text{and} 
\quad \lim_{c\to0} ~\frac{\log \alpha_{ij}^*(c)}{\log\alpha_{ks}^*(c)} = 1 .
\end{equation}
Introducing for the hypotheses~$\Hyp_i$ different observation costs~$c_i$ that may go to~$0$ at different rates, i.e., setting $A_{ji}= B_{ji}/c_i$, the results of
Theorem~\ref{Th:MSPRTnearopt} can be generalized to the more general asymmetric case where the ratios in~\eqref{symcasePrEr} are bounded away from zero and infinity. 
This generalization is important for certain applications.

Lorden's outlined results and methodologies hold significant potential for application across various problems and domains. For instance, consider their relevance in the multistream (or multichannel) problem involving two decisions and multiple data streams, as explored by \citet{FellourisTartakovsky-IEEEIT2017} and discussed in \cite[Chapter 1]{Tartakovsky_book2020}. Sequential hypothesis testing within multiple data streams, such as sensors, populations, or multichannel systems, carries numerous practical implications and applications.

Suppose observations are sequentially acquired over time in $N$ streams.
The observations  in the $i$th data stream correspond to a realization of a stochastic process $X(i)=\{X_n(i)\}_{n \in \Nbb}$, 
where $i \in \cN:=\{1, \ldots,N\}$ and  $\Nbb=\{1,2, \dots\}$.   Let  $\Hyp_{0}$ be the null  
hypothesis according to which all $N$ streams are not affected, i.e., there are no ``signals'' in all streams at all. For any given non-empty subset of components, 
$\Bc \subset \cN$, let $\Hyp_{\Bc}$ be the hypothesis according to which  only the  components $X(i)$ with $i$ in $\Bc$ contain signals. Denote by  $\Pb_{0}$ 
and $\Pb_{\Bc}$ the distributions of $\Xb$ under hypotheses $\Hyp_{0}$ and $\Hyp_{\Bc}$, respectively.
Next, let  $\ccP$ be a class of subsets of $\cN$ that incorporates {\it a priori} information that may be available regarding the subset of affected streams. Denote by 
$|\Bc|$ the size of a subset $\Bc$, i.e., the number of signals under $\Hyp_\Bc$, and by $|\ccP|$ the size of class $\ccP$, i.e., the number of possible alternatives in $\ccP$. 
For example,  if we know upper $\om \le N$ and lower $\um \ge 1$ bounds on the size of the affected subset or when we know  that at most $K$ 
streams can be affected, then   $\ccP = \ccP_{\um,\om}  =\{\Bc \subset \cN: \um \le  |\Bc| \leq \om\}$ and $\ccP = \ccP_{K}  = \{\Bc \subset \cN : 1 \le |\Bc| \le K\}$, respectively.
 
 We aim to test $\Hyp_{0}$,  the simple null hypothesis indicating no signals in any data stream, against the composite alternative
 $\Hyp_1$, according to which the subset of streams with signals belongs to  $\ccP$. We denote $\Pb_0^{n}=\Pb_0|_{\Fc_n}$ and $\Pb_\Bc^{n}=\Pb_\Bc |_{\Fc_n}$ as restrictions of probability measures $\Pb_0$ and $\Pb_\Bc$ to the $\sigma$-algebra 
$\Fc_n$, and let $p_0(\Xb^n)$ and $p_\Bc(\Xb^n)$ denote the corresponding probability densities of these measures with respect to some non-degenerate 
$\sigma$-finite measure, where $\Xb^n=(\Xb_1,\dots,\Xb_n)$ denotes the concatenation of the first $n$ observations from all data streams. 

In what follows, we confine ourselves to the i.i.d.\ scenario where observations across 
streams are independent. Moreover, within specific streams, observations are also independent, possessing densities $g_{i}(x)$ and $f_{i}(x)$ if the $i$-th stream is unaffected and contains a signal, respectively. Hence, the hypothesis testing problem can be formulated as
\begin{align*} 
\Hyp_{0} : & \quad p(\Xb^n)=p_{0}(\Xb^n) = \prod_{i=1}^N \prod_{t=1}^n g_{i}(X_t(i)); 
\\
\Hyp_{1} = \bigcup_{\Bc \in \ccP} \Hyp_{\Bc} : & \quad   p_\Bc(\Xb^n) = \prod_{i \in \Bc}\prod_{t=1}^n f_i(X_t(i)) \times \prod_{i \in \Nc \setminus \Bc}\prod_{t=1}^n g_i(X_t(i)) .
\end{align*} 
Since the hypothesis testing problem is binary the terminal decision $d$ takes two values $0$ and $1$, so $d\in \{0,1\}$ is a $\Fc_{T}$-measurable random variable 
such  that  $\{d=j\}=\{T<\infty, \Hyp_j \; \text{is selected} \}$,  $j=0,1$.  

A sequential test should be designed in such a way that the type-I (false alarm) and type-II (missed detection) 
error probabilities are controlled, i.e., do not exceed given, user-specified levels. 
Denote by  $\classPal$ the class of sequential tests with the probability of false alarm below  $\alpha_0\in(0,1)$ 
and the probability of missed detection below  $\alpha_1\in(0,1)$, i.e., 
\begin{equation} \label{ccab}
\classPal= \left\{ \delta: \Pb_{0}(d=1) \leq \alpha_0 ~~ \text{and} ~~ \max_{\Bc \in \ccP} \Pb_{\Bc}(d=0) \leq \alpha_1 \right\}.
 \end{equation}
 In general, it is not possible to design the tests that are third-order (to within $o(1)$) or even second-order (to within a constant term $O(1)$) asymptotically optimal 
 as $\alpham=\max(\alpha_0,\alpha_1)\to 0$.
 Only finding a test $T_*$ that minimizes the expected sample sizes $\Eb_0[T]$ and $\Eb_{\Bc}[T]$ for every $\Bc \in \ccP$ to first order is possible, that is,
 \begin{align*}
\Eb_{0}[T_{*}] &\sim \inf\limits_{\delta \in \classPal} \Eb_{0}[T], 
\\
\Eb_{\Bc}[T_{*}] & \sim  \inf\limits_{\delta \in \classPal} \Eb_{\Bc}[T] \quad \text{for all}  ~  \Bc \in \ccP,
\end{align*}
where $\Eb_{0}$ and $\Eb_{\Bc}$ are expectations under $\Pb_0$ and $\Pb_\Bc$, respectively.  

Hereafter we use the notation $x_\alpha \sim y_\alpha$ as $\alpha\to0$ when $\lim_{\alpha\to0} (x_\alpha/y_\alpha)=1$.  

Let $\ccP$ be an arbitrary class of subsets of $\cN$.  For any $\Bc \in \ccP$, let  $\Lambda_{\Bc}(n)$ be the likelihood ratio of $\Hyp_{\Bc}$ against $\Hyp_{0}$ 
given the observations from all streams up to time  $n$, and  let $\lambda_{\Bc}(n)$  be the corresponding log-likelihood ratio (LLR), 
\begin{align*} 
\Lambda_{\Bc} (n)  & =\frac{\drm \Pb_{\Bc}^{n}}{\drm \Pb_{0}^{n}} =  \prod_{i \in \Bc}\prod_{t=1}^n \frac{f_i(X_t(i))}{g_i(X_t(i))}  , 
\\
\lambda_{\Bc}(n)  & = \log  \Lambda_{\Bc}(n) =   \sum_{i \in \Bc} \sum_{t=1}^n \log \brcs{\frac{f_i(X_t(i))}{g_i(X_t(i))}}.
\end{align*}
The natural popular statistic for testing $\Hyp_0$ against $\Hyp_1$ at time $n$ is the maximum (generalized) likelihood ratio  (GLR) statistic 
$
 \widehat{\Lambda}(n) = \max_{\Bc \in \ccP} \;   \Lambda_{\Bc}(n). 
$
However, applying the conventional GLR statistic leads only to the first-order asymptotically optimal test. In order to obtain second and third-order optimality, we need to
modify the GLR statistic into the weighed GLR
$
 \widehat{\Lambda}(n; \pib) = \max_{\Bc \in \ccP} \;  \pi_\Bc  \Lambda_{\Bc}(n),
$
where  $\pib=\{\pi_{\Bc}, \Bc \in \ccP\}$  is a  probability mass function on $\Nc$ fully supported on $\ccP$, i.e., 
$\pi_{\Bc}>0$ for all $\Bc \in \ccP$ and $\sum_{\Bc \in \ccP} \pi_{\Bc}=1$.
The corresponding weighted generalized log-likelihood ratio (GLLR) statistic is
$
  \widehat{\lambda}(n; \pib) =\max_{\Bc \in \ccP}  \left(  \lambda_{\Bc}(n) + \log \pi_{\Bc}  \right) .
$

The \textit{Generalized Sequential Likelihood Ratio Test} (\mbox{GSLRT}) $\widehat{\delta}=(\widehat{T}, \widehat{d})$ is defined as
\begin{align*} 
  \widehat{T} =  \inf \{n \ge 1:  \widehat{\lambda}(n;\pib_{1}) \ge a_1  \; \text{or} \; \widehat{\lambda}(n;\pib_{0}) \le -a_0 \} , \quad   
  \widehat{d}=  \begin{cases}
1 & \; \text{if}  \quad \widehat{\lambda}(\widehat{T};\pib_{1}) \ge a_1   \\	
0  & \; \text{if}  \quad  \widehat{\lambda}(\widehat{T};\pib_{0}) \le -a_0 \\
\end{cases},
\end{align*}
where $\pib_j=\{\pi_{j,\Bc}, \Bc \in \cP\}$, $j=0,1$  are not necessarily identical weights and  $a_0, a_1>0$ are thresholds 
that should be selected appropriately in order to guarantee the desired error probabilities, i.e., so that  $\widehat{T}$ belongs to class $\classPal$ for given $\alpha_0$ and $\alpha_1$
with almost exact equalities. The LLR in the $i$-th stream is $\lambda_i(n) = \sum_{t=1}^n \log [f_i(X_n(i))/g_i(X_n(i))]$, so that $\lambda_\Bc(n) = \sum_{i\in \Bc} \lambda_i(n)$. 

The $\cL$-number is
\begin{equation} \label{L}
\cL_{\Bc} = \exp\left\{ - \sum_{n=1}^{\infty} \frac{1}{n} \Bigl[\Pb_{0}(\lambda_{\Bc}(n)>0)+ \Pb_{\Bc}(\lambda_{\Bc}(n) \le 0) \Bigr] \right\} ,
\end{equation}
 which takes into account the overshoot; compare with Lorden's $\cL$-numbers~\eqref{Upsilonij}.

Denote by   $\whdelta_{*}(\pib)=(\whT_*(\pib),\whd_*(\pib))$ the \mbox{GSLRT} with weights
\begin{equation} \label{choose2}
\pi_{1,\Bc}= \frac{\pi_{\Bc}}{\cL_\Bc \sum_{\Bc\in\cP}(\pi_\Bc/\cL_\Bc)} \quad \text{and} \quad  \pi_{0,\Bc}=\frac{\pi_{\Bc} \, \cL_\Bc}{\sum_{\Bc\in\cP}(\pi_\Bc \, \cL_\Bc)} , \quad \Bc \in \cP.
\end{equation}

The next theorem states that $\whdelta_{*}(\pib)$ is third-order asymptotically optimal, minimizing the weighted expected sample size $\Eb^{\pib}[T]$
to within an $o(1)$ term, where $\Eb^{\pib}$ is expectation with respect to the probability measure $\Pb^{\pib}=\sum_{\Bc\in \cP} \pi_\Bc  \, \Pb_\Bc$,
i.e., the weighted expectation $\Eb^{\pib}[\cdot] = \sum_{\Bc \in \cP} \pi_\Bc \, \Eb_\Bc[\cdot]$.

\begin{theorem} \label{Th:TOAOaverage} 
Assume the second moment conditions for LLRs $\Eb_i|\lambda_i(1)|^2 < \infty$ and  $\Eb_0|\lambda_i(1)|^2 < \infty$, $i=1,\dots,N$. Let $\alpha_0$ and $\alpha_1$ approach $0$ so that 
$|\log \alpha_0|/|\log\alpha_1| \to 1$. If thresholds $a_0$ and $a_1$ are selected so that $\whdelta_{*}(\pib)$ belongs to $\classPal$, $\Pb_0(\whd_{*}(\pib)=1) \sim \alpha_0$, and $\Pb_1(\whd_{*}(\pib)=0) \sim \alpha_1$, then the GSLRT is asymptotically optimal to third order in the class $\classPal$:
\[
\inf_{\delta \in \classPal} \Eb^{\pib}[T] =  \Eb^{\pib}[\whT_*(\pib)] +o(1) \quad \text{as}~ \alpham\to 0.
\]
\end{theorem}

The central idea of the proof of this result is to consider a purely Bayesian sequential testing problem with the $1+|\ccP|$ states ``$\Hyp_{0}: \text{density} ~ g_i$ for all  $i=1,\dots, N$'' and 
``$\Hyp_{1}^\Bc: \text{density} ~ f_\Bc$ for $\Bc\in \cP$'', and  two terminal decisions $d=0$ (accept $\Hyp_{0}$) and $d=1$ (accept $\Hyp_{1}=\cup_{\Bc\in\cP} H_{1}^\Bc$).  Then we can 
exploit Lorden's methods and results to get the proof. Without Lorden's~\citeyearpar{Lorden-AS77} paper this would not be possible. Moreover, 
the whole idea of using $\cL$-numbers for corrections is based on Lorden's fundamental contribution to the field.

\subsection{Lorden's \citeyearpar{lorden-ams70} Inequality for the Excess Over the Boundary}\label{sec:lord.ineq}

Partially motivated by seeking improved estimates of the error probabilities and other operating characteristics of Wald's SPRT discussed above, \citet{lorden-ams70} considered an upper bound for estimating a random walk's ``worst case'' expected overshoot 
\begin{equation}\label{sup.oversh}
\sup_{a\ge 0} \Eb[R_a],
\end{equation}
where $a\ge 0$ is the boundary,
\begin{equation}
R_a=S_{T(a)}- a \quad\mbox{is the overshoot,}\label{over.def}
\end{equation} 
$S_n =\sum_{t=1}^n Z_t$ is the random walk, $T(a) = \inf\{n\ge 1:\; S_n> a\}$ is the stopping time, and, relaxing slightly our notation from Section~\ref{ssec:PF}, here the $Z_n$ are i.i.d.\ random variables with positive mean~$m$; let $Z$ denote a variate with the same distribution as the $Z_n$.  Wald's \citeyearpar{wald46} equation tells us that, whenever the following quantities are finite, $m\Eb[T(a)] = \Eb[S_{T(a)}]=a+\Eb[R_a]$,
 so an upper bound on $\Eb[R_a]$ provides an upper bound on the expected stopping time~$\Eb[T(a)]$ for the random walk~$S_n$ to cross the boundary~$a$.  This is closely related to estimates of the expected stopping time~$\Eb[T_*]$ of the SPRT in \eqref{SPRT}, as we shall see below.
 
 \citet{wald47} provided the upper bound for \eqref{sup.oversh} of $\sup_{a\ge 0} \Eb[Z-a|Z>a]$, which is exact for the exponential distribution and provides reasonable bounds in some other cases, but has serious deficiencies in general: it can be difficult to calculate, is overly conservative in cases like when the distribution of $Z$ has large ``gaps,'' and may be infinite even when $\Eb[(Z^+)^2]<\infty$, a sufficient condition for finiteness of \eqref{sup.oversh}. Here and throughout this section, $z^+=\max\{z,0\}$ is the positive part of $z$.

For nonnegative $Z$, results from renewal theory \citep[see][]{feller-book_vol2_66} provide estimates of $\Eb[R_a]$ close to $\Eb[Z^2]/m=\Eb[(Z^+)^2]/m$ for both $a=0$ and as $a\rightarrow\infty$.   Lorden showed that this is indeed an upper bound for \eqref{sup.oversh}
 more generally: for arbitrary i.i.d.\ $Z_n$ allowed to be discrete or continuous, and take both positive and negative values, a necessary generalization of the renewal theory results for application to sequential testing and changepoint detection and analysis in which the $Z_n$ are log-likelihood summands or other sequential test statistic terms. 
\begin{theorem}[\citet{lorden-ams70}, Theorem~1]\label{thm:excess.iid}
If $Z,Z_1,Z_2,\ldots$ are i.i.d.\ random variables with mean~$\Eb[Z]>0$ and $\Eb[(Z^+)^2]<\infty$, then $R_a$ as defined in \eqref{over.def} satisfies
 \begin{equation}\label{lord.ineq}
\sup_{a\ge 0} \Eb[R_a]\le \frac{\Eb[(Z^+)^2]}{\Eb[Z]}.
\end{equation}
\end{theorem}

Lorden's proof of this theorem involves a number of characteristically clever  techniques, of which we highlight a few here. First, he considers the stochastic process $a\mapsto R_a$, noting that (w.p.~$1$) it is piecewise-linear, each ``piece'' having slope~$-1$. Next, since $a\mapsto R_a$ and even $a\mapsto \Eb[R_a]$ can behave erratically and be resistant to estimation and bounding, Lorden uses the smoothing technique of instead estimating $\int_0^b \Eb[R_a] \drm a$ for $b\ge 0$, which is more regularly behaved, as Lorden shows. Finally, the smoothed expected overshoot $\int_0^b \Eb[R_a] \drm a$ is bounded from above using properties of the process~$R_a$, and then bounded from below using the following sub-additivity property of the integrand $a\mapsto \Eb[R_a]$ established from the sub-additivity of $a\mapsto \Eb[T(a)]$ and  Wald's equation: For any $0\le a\le b$,
\begin{align*}
\Eb[R_a]+\Eb[R_{b-a}]&=\Eb[S_{T(a)}]-a+\Eb[S_{T(b-a)}]-(b-a)\\
&=m \Eb[T(a)]+m \Eb[T(b-a)]-b\quad\mbox{(Wald's equation)}\\
&\ge m\Eb[T(b)]-b\quad\mbox{(sub-additivity of $\Eb[T(b)]$)}\\
&=\Eb[S_{T(b)}]-b\quad\mbox{(Wald's equation)}\\
&=\Eb[R_b].
\end{align*}

Returning to the stopping time~$T_*$ of the SPRT in \eqref{SPRT}, now let the $Z_n$ be the log-likelihood ratio terms as in \eqref{llr.sprt}, $a_0$ and $a_1$ the boundaries in \eqref{SPRT}, and expectation and probability are  under the alternative hypothesis density~$f_1$. The random walk~$S_n$ now coincides with the log-likelihood ratio statistic~$\lambda_n$ in \eqref{llr.sprt}, although we continue to use the $S$ notation here for clarity. In order to relate $T_*$ to $T(a_1)$ Lorden observes that
\begin{equation}\label{ST.le.min.lord}
S_{T_*}\le\min\{S_{T_*},a_1\}+(S_{T(a_1)}-a_1),
\end{equation} and then applying \eqref{lord.ineq} to the latter term gives the upper bound
\begin{equation}\label{ET*.UB.lord}
\Eb[T_*]\le \frac{(1-\alpha_1)a_1-\alpha_1a_0}{m}+\frac{\Eb[(Z^+)^2]}{m^2}
\end{equation} on the expected stopping time of the SPRT under the alternative hypothesis, with a bound under the null hypothesis obtained analogously. \citet{wald47} provides a well-known upper bound on the type~II error probability~$\alpha_1$, but in order to apply \eqref{ET*.UB.lord} what is needed is clearly a lower bound on $\alpha_1$, and a lower bound on $\alpha_0$ for the corresponding bound under the null.  Both of these can be obtained by another application of Lorden's theorem, as follows.  Wald's argument gives that
\begin{equation*}
\frac{\alpha_0}{1-\alpha_1} = \Eb[\exp(-S_{T_*})|S_{T_*}>a_1].
\end{equation*} Using the conditional Jensen's inequality with a bound like \eqref{ST.le.min.lord} after multiplying by the indicator of the event~$\{S_{T_*}>a_1\}$, Lorden obtains
\begin{equation*}
\frac{\alpha_0}{1-\alpha_1} \ge \exp[-\Eb[(S_{T_*}|S_{T_*}>a_1)] \ge \exp\left[-\left(a_1+\frac{\Eb[(Z^+)^2]}{(1-\alpha_1)m}\right)\right].
\end{equation*} Using the standard upper bound $\alpha_1\le e^{-a_1}$, this gives
\begin{equation*}
\alpha_0  \ge (1-e^{-a_1})\exp\left[-\left(a_1+\frac{\Eb[(Z^+)^2]}{(1-e^{-a_1})m}\right)\right],
\end{equation*} with an analogous lower bound for $\alpha_1$.

\citet[][Section~2]{lorden-ams70} also obtains generalizations of \eqref{lord.ineq} to cases in which the variates~$Z_n$ are not necessarily i.i.d. They key property is the sub-additivity of $T(a)$ for which Lorden assumes the sufficient condition
\begin{equation*}
\Eb[(Z_n^+)^2|T(a)\ge n]\le r\cdot \Eb[Z_n|T(a)\ge n]
\end{equation*} for some factor~$r$. Under this condition Lorden obtains analogous bounds on $\Eb[R_a]$ and bounds on the moments $\sup_{a\ge 0} \Eb[(R_a)^p]$  for non-i.i.d.\ observations~$Z_n$ \citep[][Theorems~2 and 3]{lorden-ams70}, as well as bounds on the tail probability $\Pb(R_a>x)$ for i.i.d.\ observations \citep[][Theorem~4]{lorden-ams70}.

Other than his seminal 1971 paper on changepoint detection, \citet{lorden-ams70} is his most highly cited paper. In addition to its uses in sequential testing, changepoint detection, and renewal theory, it has found applications in reliability theory \citep{Rausand03}, clinical trial design \citep{Whitehead97}, finance \citep{Novak11}, and queuing theory \citep{Kalashnikov13}, among other applications.  Perhaps reflecting its fundamental nature and wealth of applications, Lorden's \hbox{Inequality} -- as \eqref{lord.ineq} has become known -- even has its own Wikipedia entry (\url{https://en.wikipedia.org/wiki/Lorden%27s_inequality}).


\subsection{Lorden's 2-SPRT and the Kiefer--Weiss Minimax Optimality}\label{ssec:2SPRT}

Suppose that based on a sequence of independent observations $\{X_n\}_{n \ge 1}$ with common parametric density $f_\theta$ one wishes to test
the hypothesis $\Hyp_0: \theta=\theta_0$ versus $\Hyp_1: \theta =\theta_1$ ($\theta_0< \theta_1$) with error probabilities at most $\alpha_0$ and $\alpha_1$.
Even though the SPRT has the remarkable optimality property of minimizing the expected sample size for both statistical hypotheses $\Eb_{\theta_i}[T]$, $i=0,1$,  
its performance may be poor when the true parameter value $\theta=\vartheta\in (\theta_0,\theta_1)$ differs from putative ones $\theta_0$ or $\theta_1$. 
Its expected sample size $\Eb_\vartheta[T]$ can be even much larger than that of the fixed 
sample size of the Neyman-Pearson test. See, e.g., Section 5.2 in \citet{TNB_book2014}. Much work has been directed toward finding sequential tests that reduce 
the expected sample size of the SPRT for parameter values between the hypotheses.  

Let $\class(\alpha_0,\alpha_1)=\{ \delta: \alpha_i(\delta) \leq \alpha_i, ~ i=0,1\}$ denote the class of tests with error probabilities at most $\alpha_0$ and $\alpha_1$ and let
\[
\ESS(\alpha_0,\alpha_1) =\inf_{\delta\in\class(\alpha_0,\alpha_1)} \sup_\theta \Eb_\theta[T]
\] 
denote the expected sample size of an optimal test in the class $\class(\alpha_0,\alpha_1)$ in the worst-case scenario.
The problem of finding a test~$\delta_0=(T_0,d_0)$ such that $\sup_\theta \Eb_\theta[T_0]=\ESS(\alpha_0,\alpha_1)$ subject to the error probability constraints $\alpha_0$ and $\alpha_1$
is known as the Kiefer--Weiss problem. No strictly optimal test has been found so far. \citet{Kiefer&Weiss-AMS1957} presented structured results about tests which minimize the 
expected sample size $\Eb_{\theta}[T]$ at a selected point $\theta=\vartheta\in (\theta_0,\theta_1)$, which is referred to as the modified Kiefer--Weiss problem. 
\citet{weiss-jasa62} proved that the 
Kiefer--Weiss problem reduces to the modified problem in symmetric cases for normal and binomial distributions.  \citet{lorden-as76} made a valuable contribution to the modified Kiefer--Weiss problem  for two not necessarily parametric hypotheses $\Hyp_i: \Pb=\Pb_i$, $i=0,1$, when the observations $X_1,X_2,\dots$ are i.i.d.\ and 
their true distribution $\Pb_2$ may be different from $\Pb_0$ and $\Pb_1$. \citet{lorden-as76} introduced a simple combination of one-sided SPRTs, called the
2-SPRT, and proved that it is third-order asymptotically optimal.  Later, \citet{lorden-ptrf80} proved theorems that characterize the basic structure of optimal 
sequential tests for the modified Kiefer--Weiss problem. His work has generated several works related to both the modified Kiefer--Weiss problem and 
the original Kiefer--Weiss problem of minimizing the maximal expected sample size;
see, e.g., \citet{Huffman-AS83}, \citet{DragNovikov-TVP87}, and \citet[][Section~5.3]{TNB_book2014}. 
   
Consider the following modified Kiefer--Weiss problem. Let $(\Omega, \Fc ,\Fc_n, \Pb )$, 
$n \in \Zbb_+ $, be a filtered probability space where the sub-$\sigma$-algebra $\Fc_n=\sigma(\Xb^n)$ of~$\Fc$ is 
 generated by the observations $\Xb^n = \{X_t, \: 1 \le t \le n\}$. 
The goal is to test the hypotheses~$\Hyp_i:~ \Pb = \Pb_i$, $i = 0, 1$, where $\Pb_0, \Pb_1$ are given probability measures which are locally mutually 
absolutely continuous. The true probability measure is either one of 
$\Pb_i$ or an ``intermediate'' measure $\Pb_2$ which is also locally absolute continuous with respect to $\Pb_i$. Let
$\Qb^{n}$ be a dominating measure. The observations are i.i.d.\ under $\Pb_0, \Pb_1, \Pb_2$ so the sample $\Xb^n = (X_1,\dots,X_n)$ has joint densities~$p_{i}(\Xb^n) = 
\prod_{t=1}^n f_{i}(X_t)$ for $i=0,1,2$ 
with respect to $\Qb^{n}$, where $ f_{i}(X_t)$, $t \ge 1$, are densities for the $t$-th observation.

For $n\in \Nbb$ and $i=0,1$, define the LR and LLR processes 
\begin{align*}
\Lambda_{i}(n) = \frac{\drm \Pb_2^{n}}{\drm \Pb_i^{n}}(\Xb^n) =\prod_{t=1}^n \frac{f_2(X_t)}{f_{i}(X_t)}, \quad \lambda_i(n)= \log \Lambda_i(n) = 
 \sum_{t=1}^n \log \brcs{\frac{f_2(X_t)}{f_{i}(X_t)}},
\end{align*}
with $\Lambda_i(0)=1$ and $\lambda_i(0)=0$.

Define two parallel one-sided  SPRTs
\begin{equation}\label{T01}
T_0 = \inf\set{ n\geq 1: \lambda_1(n)  \ge a_1 }, \quad T_1 = \inf\set{ n\geq 1: \lambda_0(n) \ge a_0 }  .
\end{equation}
The stopping time of Lorden's 2-SPRT~\citep{lorden-as76} is $T^\star =\min(T_0,T_1)$ and the terminal decision is $d^\star = \arg\min_{i=0,1} T_i$.   
If $a_i=\log (1/\alpha_i)$, $i=0,1$, then $\alpha_i(\delta^\star)=\Pb_i(d^\star\neq i) \leq \alpha_i$, i.e., this test belongs to class 
$\class(\alpha_0,\alpha_1)= \{ \delta: \alpha_0(\delta) \leq \alpha_0, \alpha_1(\delta) \leq \alpha_1\}$. These upper bounds may be rather conservative. For example, in the symmetric case 
$\Pb_2(d^\star=1)= \Pb_2(d^\star=0)=1/2$, we have $\alpha_i(\delta^\star)\leq \alpha_i/2$.  

Let $\Eb_{2}$ denote expectation under $\Pb_2$, and let $I_i = \Eb_2[\lambda_i(1)]$, $i=0,1$, denote Kullback--Leibler information numbers. The following theorem, proved by
\citet{lorden-as76}, establishes third-order asymptotic optimality of Lorden's 2-SPRT for small probabilities of errors $\alpha_i$. Its proof is based on Bayesian arguments.
This theorem emerges from Theorem~1 in \citet{Lorden-AS77}, which was proved a year later.

\begin{theorem}\label{Th:AOiidKW}
Let the observations $\{X_n\}_{n \geq 1}$ be i.i.d.\ under $\Pb_i$, $i=0,1,2$. Assume that the Kullback-Leibler information numbers $I_0$ and $I_1$ are positive 
and, in addition, the second-moment conditions $\Eb_{2}|\lambda_i(1)|^{2}<\infty$, $i=0,1$, hold.  Let $\alpha^\star_0(a_0,a_1)$ and $\alpha^\star_1(a_0,a_1)$ denote the
error probabilities of the 2-SPRT $\delta^\star(a_0,a_1) = (T^\star(a_0,a_1), d^\star(a_0,a_1))$. Let $\ESS(a_0,a_1)$ denote infimum of the expected sample size $\Eb_{2} [T]$
over all tests with $\Pb_0(d=1) \leq \alpha_0^\star(a_0,a_1)$ and  $\Pb_1(d=0) \leq \alpha_1^\star(a_0,a_1)$. Then
\begin{equation}\label{AO2SPRT}
 \ESS(a_0,a_1)= \Eb_{2}[T^\star(a_0,a_1)]+o(1) \quad \text{as}~ \min(a_0,a_1) \to \infty,
\end{equation}
where $o(1)\to0$ as $\min(a_0,a_1) \to \infty$.
\end{theorem}
This theorem implies that if the thresholds $a_0$ and $a_1$ in the 2-SPRT are selected so that the error probabilities $\alpha^\star_0(a_0,a_1)=\alpha_0$ 
and $\alpha^\star_1(a_0,a_1)=\alpha_1$ are exactly equal to the given values~$\alpha_0$ and $\alpha_1$, then it is third-order asymptotically optimal as $\alpham\to0$ in 
the class~$\class(\alpha_0,\alpha_1)$.  The requirement of exact error probabilities can also be relaxed to the asymptotic equalities 
$\alpha^\star_i(a_0,a_1)=\alpha_i(1+o(1))$, $i=0,1$.

\begin{figure}[ht]
\begin{center}
\includegraphics[width = 0.75\textwidth]{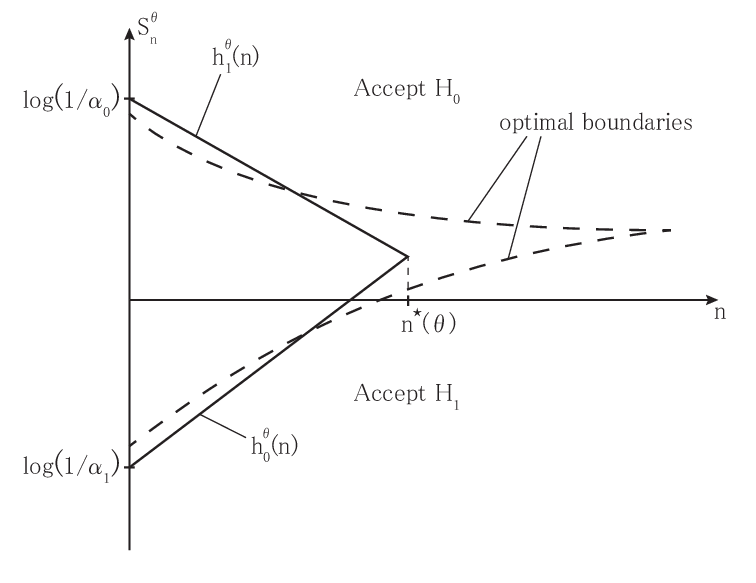}
\caption{The boundaries $h_1^\theta(n)$ and $h_0^\theta(n)$ of the $2$-SPRT (solid) and optimal boundaries (dashed) as functions of~$n$. 
$S_n^\theta= S_n - n \Eb_\theta [X_1]$. \label{fig:2SPRTtriangle}}
\end{center}
\end{figure}

The significance of this result cannot be overstated, as Lorden's simple test is nearly optimal. Simultaneously, the optimal test can be computed using Bellman's 
backward induction algorithm since the optimal sequential test is truncated, meaning it has a bounded maximal sample size, as demonstrated by \citet{Kiefer&Weiss-AMS1957}. 
For one-parameter exponential families $\{\Pb_\theta, \theta\in \theta\}$, the optimal bounds exhibit curvature in the $(S_n,n)$ plane, where $S_n =\sum_{t=1}^n X_i$, 
and determining them typically entails substantial computation. In contrast, Lorden's 2-SPRT approximates optimal curved boundaries with simple linear ones, resulting 
in a continuation region shaped like a triangle, as illustrated in Figure~\ref{fig:2SPRTtriangle}.

\citet{lorden-as76} performed an extensive performance analysis for testing the mean~$\theta$ of the Gaussian distribution $X_n \sim \Nc(\theta,1)$ 
with the hypotheses $\Hyp_i: \theta=\theta_i$, ($i=0,1$) in the symmetric case where
$\alpha_0=\alpha_1$ and $\Pb_2\sim \Nc(\theta_\star,1)$, $\theta_\star=(\theta_0+\theta_1)/2$. The conclusion is that the $2$-SPRT performs very closely to the optimal test with 
its curved boundaries 
obtained using backward induction. The efficiency depends on the error probabilities, but it was over $99\%$ in all his performed experiments. 
Similar results were obtained by \citet{Huffman-AS83} for the exponential example $f_\theta(x)=\theta e^{-\theta x}$, $x \ge 0$, $\theta >0$. 
Here the $2$-SPRT has efficiency over $98\%$ and almost always over $99\%$ for a broad range of error probabilities and parameter values.

The results of Theorem~\ref{Th:AOiidKW} can be extended to the multiple hypothesis case with $N+1$ hypotheses, $N>1$. Specifically, 
the modified matrix SPRT is also third-order asymptotically optimal as $\alpham=\max_{0\leq i\leq N}\alpha_i \to 0$ in the class of tests
\begin{equation*}
\class(\alphab)  = \set{\delta: \alpha_{i}(\delta)  \le \alpha_{i} ~\text{for} ~ i = 0,1,\dots,N},
\end{equation*}
where $\alpha_i(\delta)=\Pb_i(d \neq i)$ and $\alphab=(\alpha_{0}, \alpha_1,\dots,\alpha_N)$ is a vector of given error probabilities, c.f.\ \citet[][Theorem~5.3.3 (page 240)]{TNB_book2014}. 

Lorden's results have sparked further research into the minimax Kiefer--Weiss problem, which aims to establish near-optimal solutions for the least favorable 
intermediate distribution $\Pb_2$ within the single-parameter exponential family. Consider  the parametric case $\Pb_2=\Pb_\theta$, $\Pb_i=\Pb_{\theta_i}$ 
where the hypotheses are $\Hyp_0:\Pb_0=\Pb_{\theta_0}$ and  
$\Hyp_1:\Pb_1=\Pb_{\theta_1}$, $\theta_0 < \theta_1$. 
Let $\theta$ be an arbitrary point belonging to the interval $(\theta_0,\theta_1)$ and let $\delta^\star(\theta)=(d^\star(\theta),T^\star(\theta))$ denote the $2$-SPRT tuned to~$\theta$. 
In other words, $T^\star(\theta)=\min(T_0^\theta,T_1^\theta)$, where the $T_i^\theta$'s are defined by~\eqref{T01} with the LLRs 
$\lambda_i^\theta(n)= \log [\drm\Pb_\theta^n/\drm\Pb_{\theta_i}^n](\Xb^n)$, $i=0,1$, tuned to~$\theta$. 
 
 Theorem~\ref{Th:AOiidKW} implies that the $2$-SPRT $\delta^\star(\theta)$ is third-order asymptotically optimal for minimizing $\Eb_\theta [T]$ at the intermediate point 
 $\theta\in(\theta_0,\theta_1)$ when the second moments $\Eb_\theta|\lambda_0^\theta(1)|^2$ and $\Eb_\theta|\lambda_1^\theta(1)|^2$ are finite, and the thresholds $a_i$ are selected so that the error probabilities are either exactly equal, or at least close, to the given numbers $\alpha_i$, a challenging task. However, setting $a_i=|\log\alpha_i|$ embeds the $2$-SPRT into the class $\class(\alpha_0,\alpha_1)$, and Theorem~\ref{Th:AOiidKW} suggests that if one can find a nearly least favorable point $\theta^*$, i.e., $\theta^*$ can be selected so that $\sup_\theta \Eb_\theta[T^\star(\theta^*)] \approx \Eb_{\theta^*}[T^\star(\theta^*)]$, then $\delta^\star(\theta^*)$ is an approximate solution to the Kiefer--Weiss problem of minimizing $\sup_\theta \Eb_\theta[T]$.

For the single-parameter exponential family with density $\{f_\theta(x), \theta\in\Theta\}$, where
\begin{equation}\label{Expfam}
\frac{f_\theta(x) }{f_{\tilde\theta}(x)} = \exp\set{(\theta-\tilde\theta) x -(b(\theta) - b(\tilde\theta))},
\end{equation}
with~$b(\theta)$ being a convex and infinitely differentiable function on~$\widetilde\Theta\subset \Theta$,
it is feasible to identify the nearly least favorable point $\theta^*(\alpha_0,\alpha_1,\theta_0,\theta_1)$ such that the 
$2$-SPRT with thresholds $a_i=\log(1/\alpha_i)$ achieves second-order asymptotic minimaxity, meaning the residual term in the discrepancy between the expectation 
of the sample size of the optimal test and the \mbox{2-SPRT} is of order $O(1)$ for small $\alpha_i$. Initially addressed by \citet{Huffman-AS83}, who proposed 
$\theta^*(\alpha_0,\alpha_1,\theta_0,\theta_1)$ leading to a residual term of order $o(|\log \alpha|^{1/2})$, 
this problem was further advanced by \citet{DragNovikov-TVP87}, who demonstrated the second-order optimality of Huffman's version of the 2-SPRT.

As previously noted, the formulas $a_i=|\log\alpha_i|$, ensuring the inequalities $\alpha_i(\delta^\star(\theta)) \le \alpha_i$, tend to be overly conservative. A refinement 
can be achieved by observing that
\begin{align*}
\alpha_1(\delta^\star(\theta)) & = \Pb_\theta(T^\star =T_0^\theta) e^{-a_1} \Eb_\theta\set{e^{-(a_1-\lambda_1^\theta(T_0^\theta))} | T^\star=T_0^\theta } ,
\\
\alpha_0(\delta^\star(\theta)) & = \Pb_\theta(T^\star =T_1^\theta) e^{-a_0} \Eb_\theta\set{e^{-(a_0-\lambda_0^\theta(T_1^\theta))} | T^\star=T_1^\theta }
\end{align*}
where, asymptotically as $a_i \to \infty$,
\[
  \Eb_\theta\set{e^{-(a_0-\lambda_0^\theta(T_1^\theta))} | T^\star=T_1^\theta } \to \zeta_0^\theta, \quad
  \Eb_\theta\set{e^{-(a_1-\lambda_1^\theta(T_0^\theta))} | T^\star=T_0^\theta } \to \zeta_1^\theta.
\]
For the non-arithmetic case, $\zeta_i^\theta$ can be computed using the renewal-theoretic argument similar to \eqref{Upsilonij} and \eqref{zetaell}:
\[
\zeta_{i}^\theta = \frac{1}{I(\theta,\theta_i)}  \exp\set{- \sum_{n=1}^\infty \frac{1}{n}  \brcs{\Pb_\theta (\lambda_{i}^\theta(n) >0)  + \Pb_i(\lambda_{i}^\theta(n) \le 0)}} ,
\]
where 
$
I(\theta, \theta_i) = (\theta-\theta_i) \dotb{b}(\theta)  -(b(\theta) - b(\theta_i)) 
$
are the Kullback--Leibler numbers.

The efficiency of Lorden's 2-SPRT in stopping to both reject and accept the null was appealing in clinical trial designs, where these actions are known as efficacy and futility stopping. \citet{Lai&Shih} extended the 2-SPRT from the fully sequential setting to the group sequential setting while maintaining its efficiency, and balancing the tradeoff between efficacy and futility stopping.

\subsection{Near Uniform Optimality of the GLR SPRT for Composite Hypotheses} \label{ssec:GLR}

For practical purposes, it is considerably more significant to devise tests that minimize the expected sample size $\Eb_\theta[T]$  for all possible parameter values 
(i.e., uniformly optimal) rather than to address the minimax Kiefer--Weiss problem of minimizing $\Eb_\theta[T]$ at a least favorable point. In this section, our primary objective is to explore the design of sequential tests that are at least approximately uniformly optimal for small error probabilities or asymptotically Bayesian for a small cost of observations for testing composite hypotheses.

Consider a sequence of i.i.d.\ observations  $X_1,X_2,\dots$ originating from a common distribution $\Pb_\theta$ with density~$f_\theta(x)$ with respect to some 
non-degenerate sigma-finite measure, where the $\ell$-dimensional parameter $\theta=(\theta_1,\dots,\theta_\ell)$ 
belongs to a subset $\Theta$ of the Euclidean space~$\Rbb^\ell$. The parameter space~$\Theta$ is partitioned into $3$ disjoint sets $\Theta_0, \Theta_1$ and $\Iin$, i.e., 
$\Theta=\Theta_0 \cup \Theta_2 \cup \Iin$. The objective is to test the two composite hypotheses~$\Hyp_0: \theta \in \Theta_0$ against $\Hyp_1: \theta \in \Theta_1$.
The subset~$\Iin$ of~$\Theta$ denotes an indifference zone where the loss~$L(\theta, d)$ associated with correct or incorrect decisions~$d$ is zero, i.e., 
no constraints on the probabilities~$\Pb_\theta(d=i)$ are imposed if $\theta\in \Iin$. The introduction of an indifference zone is typically motivated by the recognition that in many applications, the correct action is not crucial and often not even feasible when the hypotheses are very close. However, in principle~$\Iin$ may be an empty set.

We aim to find a sequential test $\delta=(T,d)$ that minimizes the expected sample size~$\Eb_\theta [T]$ uniformly for all $\theta \in \Theta$ in the class 
of tests~$\class(\alpha_0, \alpha_1)$  in which the maximal error probabilities $\sup_{\theta \in \Theta_i} \Pb_\theta(d \neq i)$ 
are upper-bounded by the given values:  
\begin{equation}\label{classcomp2hyp}
\class(\alpha_0,\alpha_1) = \set{\delta: \sup_{\theta\in \Theta_0} \Pb_\theta(d=1) \le \alpha_0 ~~ \text{and} ~~ \sup_{\theta\in \Theta_1} \Pb_\theta(d=0) \le \alpha_1}.
\end{equation}
Thus, we are interested in the frequentist problem of finding a test $\delta_{\rm opt}$ such that
\begin{equation}\label{OptuniformProblem}
\inf_{\delta\in\class(\alpha_0,\alpha_1)} \Eb_\theta [T] = \Eb_\theta [T_{\rm opt}]\quad \text{uniformly in}~~ \theta\in \Theta.
\end{equation}
Unfortunately, such a uniformly optimal solution does not exist, and one has to resort to finding asymptotic approximations for small error probabilities. In the frequentist setting, it 
is possible to find first-order asymptotically optimal tests that satisfy
\begin{equation} \label{FOopt}
\lim_{\alpham \to 0} \frac{\inf_{\delta\in\class(\alpha_0,\alpha_1)} \Eb_\theta [T]}{\Eb_\theta [T]} =1 \quad \text{for all}~~ \theta\in\Theta.
\end{equation}

In addition to  the frequentist problems~\eqref{OptuniformProblem}-\eqref{FOopt}, it is of interest  to consider a Bayesian approach putting an {\it a priori} 
distribution~$W(\theta)$ on~$\Theta$ with a cost~$c$ per observation and a 
loss function~$L(\theta)$ at the point~$\theta$ associated with accepting the incorrect hypothesis and find asymptotically optimal
tests when the cost $c$ is small. The Bayes average (integrated) risk of a sequential test $\delta=(T,d)$ is
\[
\rho_c^W(\delta) = \int_{\theta\le \theta_0} L(\theta) \Pb_\theta(d=1) \, W(\drm\theta) + \int_{\theta \ge \theta_1} L(\theta) \Pb_\theta(d=0) \, W(\drm\theta) +
c \int_{\Theta} \Eb_\theta [T] \, W(\drm\theta) .
\]
It turns out that in the Bayesian context, it is possible to find tests that are not only asymptotically (as $c\to0$) first-order optimal, 
$\inf_\delta \rho_c^W(\delta)= \rho_c^W(\delta)(1+o(1))$, but also second-order optimal, i.e., $\inf_\delta \rho_c^W(\delta)= \rho_c^W(\delta) + O(c)$ and even third-order optimal, i.e., 
$\inf_\delta \rho_c^W(\delta)= \rho_c^W(\delta) + o(c)$.

In the case of the one-parameter exponential family~\eqref{Expfam}, using optimal stopping theory, it can be shown that the optimal Bayesian test 
$\delta_{\rm opt}=(T_{\rm opt}, d_{\rm opt})$ is
\[
T_{\rm opt} = \inf\set{n\ge 1: (S_n,n) \in \cB_c}, \quad d_{\rm opt}= j~~\text{if}~(S_n,n) \in \cB_c^j, \quad  j =0,1,
\]
where $S_n=X_1+\cdots+X_n$ and $\cB_c =\cB_c^0 \cup\cB_c^1$ is a set that can be found numerically.

\citet{Schwarz-AMS1962} derived the test~$\delta^\star(\hat\theta)$ with $\hat\theta=\{\hat\theta_n\}$ being the maximum likelihood estimator (MLE) of~$\theta$, 
as an asymptotic solution as $c\to 0$ to the Bayesian problem with the $0-1$~loss function.
Specifically, the  {\em a posteriori} risk of stopping is
\begin{equation}\label{APRcomp}
R_n^{\rm st}(S_n) =\min_{i=0,1} \set{\frac{\int_{\Theta_i} \exp\set{\theta S_n -n b(\theta)} \, W(\drm\theta)}{\int_{\Theta} \exp\set{\theta S_n -n b(\theta)} \, W(\drm\theta)}} ,
\end{equation}
 where $\Theta_0=\{\theta \le \theta_0\}$, $\Theta_1=\{\theta \ge \theta_1\}$.  Schwarz showed that $\cB_c/|\log c| \to \cB_0$ as~$c\to 0$ and proposed a simple procedure: 
 continue sampling until $R_n^{\rm st}(S_n)$ is less than~$c$ and upon stopping accept the hypothesis for which the minimum is attained in~\eqref{APRcomp}. 
 Denote this procedure by $\widetilde\delta(c)=(\widetilde{T}(c), \widetilde{d}(c))$.
Applying Laplace's asymptotic integration method to evaluate the integrals in~\eqref{APRcomp} leads to the likelihood ratio test where the true parameter is replaced by 
the MLE~$\hat\theta_n$. This approximation prescribes stopping sampling at the time 
$\widehat{T}(\hat\theta)= \min(\widehat{T}_0(\hat\theta),\widehat{T}_1(\hat\theta))$, where
\begin{equation}\label{hatTiSchwarz}
\begin{aligned}
\widehat{T}_i(\hat\theta) 
= \inf\set{n: \sup_{\theta\in\Theta} [\theta S_n - n b(\theta)] -  [\theta_i S_n - n b(\theta_i)] \ge | \log c|} .
\end{aligned}
\end{equation}
The terminal decision rule~$\hat{d}(\hat\theta)$ of the test 
$\hat\delta(\hat\theta)= (\widehat{T}(\hat\theta), \hat{d}(\hat\theta))$ accepts~$\Hyp_0$ if $\hat\theta_{\widehat{T}} < \theta^*$, where~$\theta^*$ is such that 
$I(\theta^*,\theta_0) = I(\theta^*, \theta_1)$. Note also that
\begin{equation}\label{SchST}
\widehat{T} = \inf\set{n\ge 1 :  n \max[I(\hat\theta_n,\theta_0), I(\hat\theta_n,\theta_1)] \ge |\log c|}.
\end{equation}

The tests which use the maximum likelihood estimators of unknown parameters are usually referred to as the {\em Generalized Sequential Likelihood Ratio Tests} (GSLRT). 

\citet{WongAMS-1968} showed that the GSLRT $\hat\delta$ is first-order asymptotically Bayes as~$c \to 0$:
\[
\rho_c^W(\hat{\delta}) \sim \inf_\delta \rho_c^W(\delta) \sim c|\log c| \int_\Theta \frac{W(\drm \theta)}{I_{\max}(\theta)} , \quad 
\Eb_\theta [\widehat{T}] \sim \frac{|\log c|}{I_{\max}(\theta)} \quad \text{for every}~~ \theta \in \Theta,
\]
where $I_{\max}(\theta)=\max\set{I(\theta,\theta_0), I(\theta,\theta_1)}$.

\citet{KieferSacks-AMS1963} showed that the procedure $\widetilde\delta(c)=(\widetilde{T}(c), \widetilde{d}(c))$ with the stopping time
$
\widetilde{T}(c) = \inf\set{n \ge 1:  R_n^{\rm st}(S_n) \le c},
$
proposed by \citet{Schwarz-AMS1962}, is also first-order asymptotically Bayes. In other words, for any prior distribution~$W$,  $\rho_c^W(\widetilde\delta(c))$ behaves asymptotically like $ \inf_\delta \rho_c^W(\delta)$
as~$c \to 0$.   \citet{Lorden1} refined this result by introducing the stopping region as the first~$n$ such that $R_n^{\rm st}(S_n) \leq Q c$, where $Q$ is a positive constant, 
and demonstrated that it can be made second-order asymptotically optimal, i.e., $\inf_\delta \rho_c^W(\delta)=\rho_c^W(\widetilde\delta(Qc)) +O(c)$ as~$c\to 0$, 
while $\inf_\delta \rho_c^W(\delta)=O(c|\log c|)$. It's noteworthy that the problem addressed by  \citet{Lorden1} is more general than what we are discussing here since it encompasses general i.i.d.\ models, not limited to exponential families, and multiple-decision cases. 
Additionally, see \citet{Lorden-ams72} for multiple hypotheses in one-parameter exponential families.

A significant advancement in Bayesian theory for testing separated hypotheses about the parameter of the one-parameter exponential family \eqref{Expfam} was made by \citet{Lorden-unpublished-1977} (an unpublished manuscript). In this work, Lorden demonstrated that the family of GSLRTs can be devised to ensure third-order asymptotic optimality. This implies that they achieve the Bayes risk to within $o(c)$ as~$c\to0$.

Lorden provided sufficient conditions for families of tests to be third-order asymptotically Bayes and presented examples of such procedures based not only on the Generalized Likelihood Ratio (GLR) approach but also on mixtures of likelihood ratios. Furthermore, the error probabilities of the GSLRTs were evaluated asymptotically as a consequence of a general theorem on boundary-crossing probabilities.

Due to the significance of this work, let's delve into a more detailed overview of Lorden's theory. It's worth noting that the paper by \citet{Lorden-unpublished-1977} extends the results obtained by \citet{Lorden-AS77} for multiple discrete cases, which we discussed in Subsection~\ref{ssec:MSPRT}, to the continuous parameter case.

The hypotheses to be tested are $\Hyp_0: \underline{\theta} \leq \theta \leq \theta_0$ and $\Hyp_1: \overline{\theta} \geq \theta \geq \theta_1$, where $ \underline{\theta}$ 
and~$\overline{\theta}$ are interior points of the natural parameter space~$\Theta$. Let $\hat\theta_n\in [ \underline{\theta}, \overline{\theta} ]$ be the MLE that maximizes the 
likelihood over~$\theta$ in~$ [ \underline{\theta}, \overline{\theta} ]$. Lorden's GSLRT stops at~$\widehat{T}$  which is the minimum of the Markov times 
$\widehat{T}_0, \widehat{T}_1$ defined as
\begin{equation}\label{hatTiLorden}
\begin{aligned}
\widehat{T}_0(\hat\theta) & = \inf\set{n \geq 1: \sum_{k=1}^n \log \brcs{\frac{f_{\hat\theta_n}(X_k)}{f_{\theta_0}(X_k)} ~h_0(\hat\theta_n)} \geq a ~\text{and}~\hat\theta_n \geq \theta^*} ,
\\
\widehat{T}_1(\hat\theta) & = \inf\set{n \geq 1: \sum_{k=1}^n \log \brcs{\frac{f_{\hat\theta_n}(X_k)}{f_{\theta_1}(X_k)} ~h_1(\hat\theta_n)} \geq a ~\text{and}~\hat\theta_n \leq \theta^*} ,
\end{aligned}
\end{equation}
where~$a$ is a threshold,  $\theta^*$ satisfies $I(\theta^*,\theta_0)=I(\theta^*,\theta_1)$, and $h_0, h_1$ are positive continuous functions on $[\theta^*, \overline{\theta}], [\underline{\theta}, \theta^*]$, respectively. The hypothesis~$\Hyp_i$ is rejected when $\widehat{T}=\widehat{T}_i$.   To summarize, Lorden's family of GSPRTs is defined as
\begin{equation}\label{GSLRTlorden}
\widehat{T}(\hat\theta) = \min \set{\widehat{T}_0(\hat\theta), \widehat{T}_1(\hat\theta)}, \quad
\hat{d}= \begin{cases}
0 & \text{if}~~ \widehat{T}(\hat\theta)=\widehat{T}_1(\hat\theta)
\\
 1& \text{if}~ \widehat{T}(\hat\theta)=\widehat{T}_0(\hat\theta)
 \end{cases} ,
\end{equation}
with the $\widehat{T}_i(\hat\theta)$'s in~\eqref{hatTiLorden}.

Denote by
\[
\lambda_n(\theta,\theta_i) =  \sum_{k=1}^n \log \brcs{\frac{f_{\theta}(X_k)}{f_{\theta_i}(X_k)}} = (\theta-\theta_i) S_n - [b(\theta) - b(\theta_i)] n
\]
the LLR between points $\theta$ and $\theta_i$.

Lorden assumes that the prior distribution $W(\theta)$ has a continuous density~$w(\theta)$ positive on $[ \underline{\theta}, \overline{\theta} ]$, and that the loss~$L(\theta)$ equals zero in the indifference zone $(\theta_0,\theta_1)$ and is continuous and positive elsewhere and bounded away from~$0$ on $[\underline{\theta}, \theta_0]\cup [\theta_1, \overline{\theta}]$.
The main results in \cite[Theorem~1]{Lorden-unpublished-1977} can be briefly outlined as follows.
\begin{description}
\item [(i)] Under these assumptions the family of GSLRTs defined by \eqref{hatTiLorden}--\eqref{GSLRTlorden} with $a= |\log c| -\tfrac{1}{2} \log |\log c|$ is second-order asymptotically optimal, i.e.,
\[
\rho_c^w(\hat\delta) = \inf_{\delta} \rho_c^w(\delta) + O(c) \quad \text{as}~ c \to 0,
\]
where $\rho_c^w(\delta)$ is the average risk of the test~$\delta=(T,d)$:
\[
\rho_c^w(\delta) = \int_{\underline{\theta}}^{\theta_0} L(\theta) \Pb_\theta(d=1) w(\theta) \, \drm\theta + \int_{\theta_1}^{\overline{\theta}} L(\theta) \Pb_\theta(d=0) w(\theta) \, \drm\theta +
c \int_{\underline{\theta}}^{\overline{\theta}} \Eb_\theta [T] w(\theta)\, \drm\theta .
\]

\item [(ii)] This result can be improved from~$O(c)$ to~$o(c)$, i.e., to the third order
\[
\rho_c^w(\hat\delta) = \inf_{\delta} \rho_c^w(\delta) + o(c) \quad \text{as}~ c \to 0,
\]
making the right choice of the functions~$h_0$ and~$h_1$ by setting
\[
h_i(\theta) = \sqrt{\frac{2\pi}{I^3(\theta,\theta_i) \ddot{b}(\theta)}} ~\frac{w(\theta)|\dotb{b}(\theta)-\dotb{b}(\theta_i)|}{w(\theta_i) L(\theta_i)\zeta(\theta,\theta_i)}, \quad i=0,1 ,
\]
where $\zeta(\theta,\theta_i)=\Lc(\theta,\theta_i)/I(\theta,\theta_i)$ is a correction for the overshoot over the boundary, the factor which is the subject of renewal theory.
Specifically,
\begin{equation}\label{zetatheta}
\zeta(\theta,\theta_i)= \lim_{a\to\infty} \Eb_\theta \exp\set{-[\lambda_{\tau_a}(\theta,\theta_i) -a]}, \quad \tau_a =\inf\set{n: \lambda_n(\theta,\theta_i) \geq a},
\end{equation}
where in the non-arithmetic case $\zeta(\theta,\theta_i)$ can be computed as
\begin{equation}\label{Upsilontheta}
\zeta(\theta,\theta_i) = \frac{1}{I(\theta,\theta_i)} \exp\set{- \sum_{n=1}^\infty \frac{1}{n} \brcs{\Pb_\theta(\lambda_n(\theta,\theta_i) \leq 0) + \Pb_{\theta_i}(\lambda_n(\theta,\theta_i) > 0) }}.
\end{equation}
 Since the Bayes average risk $\inf_\delta \rho_c^w(\delta)$ is of order~$c|\log c|$, this implies that the asymptotic relative efficiency  
 $\cE_c= [\rho_c^w(\hat\delta)-\inf_\delta \rho_c^w(\delta)]/\rho_c^w(\hat\delta)$ of Lorden's test is of order $1-o(1/|\log c|)$ as~$c\to 0$.
\end{description}

Note the crucial difference between Schwarz's GSLRT~\eqref{hatTiSchwarz} and Lorden's GSLRT \eqref{GSLRTlorden}. In the Schwarz test, $h_i\equiv 1$ and the threshold 
is set as $a=|\log c|$. However, in the Lorden test, two innovations emerge. Firstly, the threshold is reduced by $\tfrac{1}{2} \log |\log c|$, and secondly,  adaptive weights~$h_i(\hat\theta_n)$ are incorporated into the GLR statistic. Since the stopping times~$\widehat{T}_i$ can be written as
\begin{equation}\label{hatTi2}
\begin{aligned}
\widehat{T}_0(\hat\theta) & = \inf\set{n \geq 1: \lambda_n(\hat\theta_n,\theta_0)  \geq a -\log h_0(\hat\theta_n)~\text{and}~\hat\theta_n \geq \theta^*} ,
\\
\widehat{T}_1(\hat\theta) & = \inf\set{n \geq 1: \lambda_n(\hat\theta_n,\theta_1)  \geq a -\log h_1(\hat\theta_n)~\text{and}~\hat\theta_n \leq \theta^*} ,
\end{aligned}
\end{equation}
Lorden's GSLRT can alternatively be perceived as the GSLRT with curved adaptive boundaries
\[
a_i(\hat\theta_n) =  |\log c| -\tfrac{1}{2} \log |\log c| - \log h_i(\hat\theta_n), \quad i=0,1,
\]
which depend on the behavior of the MLE~$\hat\theta_n$. These two innovations render this modification of the GLR test nearly optimal.

Given the complexity of Lorden's formal mathematical proof, we offer a heuristic sketch that captures the main ideas of the approach. The Bayesian perspective naturally guides us toward the mixture LR statistics
\[
\bar{\Lambda}_n^i = \frac{\int_{\underline{\theta}}^{ \overline{\theta}} e^{\theta S_n- n b(\theta)} w(\theta) \, \drm \theta}{\int_{\Theta_i} L(\theta) e^{\theta S_n- n b(\theta)} w(\theta)\, \drm \theta} , \quad i=0,1,
\]
where $\Theta_0= [\underline{\theta}, \theta_0]$, $\Theta_1=[\theta_1, \overline{\theta}]$ and $L(\theta)=1$ for the simple $0-1$~loss function. Indeed, the {\mbox{{\em a posteriori} } stopping risk is given by
\begin{equation}\label{APRcomp2}
R_n^{\rm st}(S_n) =\min_{i=0,1} \set{\frac{\int_{\Theta_i} L(\theta) e^{\theta S_n -n b(\theta)} w(\theta)  \, \drm \theta}{\int_{\underline{\theta}}^{ \overline{\theta}} e^{\theta S_n -n b(\theta)}  w(\theta)  \, \drm \theta}} .
\end{equation}
A candidate for the approximate optimum is the procedure that stops as soon as $R_n^{\rm st}(S_n) \leq A_c$ for some~$A_c \approx c$. This is equivalent to stopping as soon as $\max_{i=0,1} \bar\Lambda_n^i \geq 1/A_c$. The GLR statistics are approximated as
\[
\hat\Lambda_n^i = \frac{\max_{\theta\in [\underline{\theta}, \overline{\theta}]} ~e^{\theta S_n- n b(\theta)}}{\max_{\theta\in\Theta_i} ~e^{\theta S_n- n b(\theta)}} \approx  \frac{\max_{\theta\in [\underline{\theta}, \overline{\theta}]} ~e^{\theta S_n- n b(\theta)}}{e^{\theta_i S_n- n b(\theta_i)}} , \quad i=0,1,
\]
and the stopping posterior risk~\eqref{APRcomp2} is approximated as
\begin{equation}\label{APRcompapprox}
R_n^{\rm st}(S_n) \approx \min_{i=0,1} \frac{w(\theta_i) L(\theta_i) [\ddot{b}(\hat\theta_n)/2\pi n]^{1/2}}{w(\hat\theta_n) |\dotb{b}(\theta_i) - \dotb{b}(\hat\theta_n)|} ~e^{-\lambda_n(\hat\theta_n,\theta_i)},
\end{equation}
where $i=0$ if $\hat\theta_n \leq \theta^*$ and $i=1$ otherwise. These approximations stem from Laplace's method for asymptotic integral expansions, and its variations.

Subsequently, Lorden demonstrated the existence of $Q>1$ such that if the stopping risk exceeds $Q c$, then the continuation risk becomes smaller than the stopping risk. 
Therefore, it is approximately optimal to stop at the first instance such that $R_n^{\rm st}$ falls below~$Q c$. This finding, coupled with the approximation \eqref{APRcompapprox}, results in $T_{\rm opt}\approx\min(\tau_0, \tau_1)$, where
\[
\tau_i=  \inf\set{n: e^{-\lambda_n(\hat\theta_n,\theta_i)}/ \widetilde{h}_i(\hat\theta_n)n^{1/2} \leq  Qc} = \inf\set{n: \lambda_n(\hat\theta_n,\theta_i)\geq -\log[n^{1/2} \widetilde{h}_i(\hat\theta_n)  Qc]}
\]
with $\widetilde{h}_i(\hat\theta_n)$ given by
\[
\widetilde{h}_i(\hat\theta_n) = \sqrt{\frac{2\pi}{\ddot{b}(\hat\theta_n)}}\frac{w(\theta)|\dotb{b}(\hat\theta_n)-\dotb{b}(\theta_i)|}{w(\theta_i) L(\theta_i)} .
\]
For small~$c$, the expectation $\Eb_\theta [\tau_i]$ is of order~$|\log c|$, so $n^{1/2}$ can be replaced by~$|\log c|^{1/2}$, which yields
\[
\tau_i \approx \widehat{T}_i =  \inf\set{n: \lambda_n(\hat\theta_n,\theta_i)\geq - \log[c |\log c|^{1/2} Q \widetilde{h}_i(\hat\theta_n)]} .
\]
Note that these stopping times look exactly like the ones defined in~\eqref{hatTi2} with the stopping boundaries
\[
a_i(\hat\theta_n) =  |\log c| -\tfrac{1}{2} \log |\log c| - \log [Q \widetilde{h}_i(\hat\theta_n)], \quad i=0,1.
\]
The test based on these stopping times is already optimal to the second order. However, to achieve third-order optimality, one must carefully choose the constant $Q$ to address the overshoots. Specifically, leveraging this result, Lorden demonstrates that the risks of an optimal rule and of the GSLRT are both linked to the risks of the family of one-sided tests $\tau_a(\theta) =\inf\{n: \lambda_n(\theta,\theta_i) \geq a\}$, which are strictly optimal in the problem:
\[
\rho(\theta, v) = \inf_{T} \set{\Eb_\theta T + v \Pb_{\theta_i}(T<\infty)} =  \inf_{T} \Eb_\theta \set{ T + v \prod_{n=1}^T\frac{p_{\theta_i}(X_n)}{p_\theta(X_n)}}.
\]
If we set $a= \log[v \Lc(\theta,\theta_i)]$, then by taking $Q=1/\Lc(\theta,\theta_i)$ the resulting test will be nearly optimal to within~$o(c)$. 
Since $\theta$ is unknown, we need to replace it with the estimate~$\hat\theta_n$ to obtain
\[
a_i(\hat\theta_n) =  |\log c| -\tfrac{1}{2} \log |\log c| - \log [\widetilde{h}_i(\hat\theta_n)/\Lc(\hat\theta_n,\theta_i)] = |\log c| -\tfrac{1}{2} \log |\log c| - \log [h_i(\hat\theta_n)].
\]

It's intriguing to compare Lorden's approach with the Kiefer--Sacks test that stops the first time $R_n^{\rm st}$ becomes smaller than~$c$. 
Lorden's approach allows us to show that the test with the stopping time
\[
\widehat{T} = \inf\set{n: R_n^{\rm st}(S_n) \leq  c/\Lc(\hat\theta_n)},
\]
where $\Lc(\hat\theta_n) = \Lc(\hat\theta_n,\theta_1)$ if $\hat\theta_n < \theta^*$ and $\Lc(\hat\theta_n)=\Lc(\hat\theta_n,\theta_0)$ otherwise, is nearly optimal 
to within~$o(c)$. It's worth recalling that the factor $\zeta(\theta,\theta_i)=I(\theta,\theta_i)^{-1}\Lc(\theta,\theta_i)$ provides a necessary correction for the excess over the thresholds at stopping; see \eqref{zetatheta}.
This offers a significant enhancement over the Kiefer--Sacks test, which disregards the overshoots. Notably, this improvement is not limited to testing close hypotheses when $\Lc(\theta,\theta_i) \ll 1$. Even in cases where the parameter values are well-separated, this correction could be crucial.  For instance, in the binomial case with the success 
probabilities~$\theta_1=0.6$ and~$\theta_0=0.4$, we have $\Lc(\theta_1,\theta_0) \approx 1/15$, so Lorden's test will terminate much earlier.

Certainly, it's important to note that implementing Lorden's fully optimized GSLRT may encounter difficulties. This is primarily because computing the numbers $\zeta(\theta,\theta_i)$ analytically is often not feasible, except for specific models such as the exponential. For instance, when testing the mean in the Gaussian case, these numbers can only be computed numerically. While Siegmund's~(\citeyear{siegmund-book85}) corrected Brownian motion approximations can be utilized, they are sufficiently accurate only when the difference between $\theta$ and $\theta_i$ is relatively small. Hence, for practical purposes, only partially optimized solutions, which provide $O(c)$-optimality, are typically feasible. A workaround involves discretizing the parameter space.

Let $\hat\alpha_0(\theta)=\Pb_\theta(\hat{d}=1)$, $\theta\in \Theta_0= [\underline{\theta}, \theta_0]$ and $\hat\alpha_1(\theta)=\Pb_\theta(\hat{d}=0)$, $\theta\in \Theta_1=[\theta_1, \overline{\theta}]$ denote
the error probabilities of the GSLRT~$\hat\delta_a$.  Note that due to the monotonicity of~$\hat\alpha_i(\theta)$, $\sup_{\theta\in\Theta_i} \hat\alpha_i(\theta) = \hat\alpha_i(\theta_i)$.
In addition to the Bayesian third-order optimality property, Lorden established asymptotic approximations to the error probabilities of the GSLRT. Specifically, by Theorem~2 of
\citet{Lorden-unpublished-1977},
\begin{equation}\label{PEGSLRTlorden}
\hat\alpha_i(\theta_i)  = \sqrt{a} e^{-a} C_i(\theta_i) (1+o(1)), \quad i=0,1 \quad \text{as}~a \to \infty ,
\end{equation}
where
\[
\begin{aligned}
C_0(\theta_0) &= \int_{\theta^*}^{\overline{\theta}} \zeta(\theta,\theta_0) h_0(\theta) \sqrt{\frac{\ddot{b}(\theta)}{2\pi I(\theta,\theta_0)}} \, \drm\theta ,
\\
C_1(\theta_1) &= \int_{\underline{\theta}}^{\theta^*} \zeta(\theta,\theta_1) h_1(\theta) \sqrt{\frac{\ddot{b}(\theta)}{2\pi I(\theta,\theta_1)}} \, \drm\theta
\end{aligned}
\]
and where $\zeta(\theta,\theta_i)$, $i=0,1$, are defined in \eqref{zetatheta}--\eqref{Upsilontheta}. These approximations hold significance for frequentist problems, which are typically of primary interest in most applications. While there are no strict upper bounds on the error probabilities, leading to no specific prescription on how to embed the 
GSLRT into class~$\class(\alpha_0,\alpha_1)$, the asymptotic approximations~\eqref{PEGSLRTlorden} enable us to select thresholds~$a_i$ in the stopping times~$\widehat{T}_i$ 
so that $\hat\alpha_i(\theta_i)\approx \alpha_i$, $i=0,1$, at least for sufficiently small~$\alpha_i$.  Note that in this latter case, the threshold 
$a$ in~\eqref{hatTi2} should be replaced with~$a_i$, the roots of the transcendental equations
\[
a_i - \frac{1}{2} \log a_i  = \log [C_i(\theta_i)/\alpha_i], \quad i=0,1.
\]
With this choice, the GSLRT is asymptotically uniformly  first-order optimal with respect to the expected sample size, i.e.,
\[
\inf_{\delta\in \class(\alpha_0,\alpha_1)} \Eb_\theta [T] = \Eb_\theta [\widehat{T}] (1+o(1)) \quad \text{as}~ \alpha_{\max} \to 0 \quad 
\text{for all}~\theta \in [\underline{\theta}, \overline{\theta}],
\]
where the $o(1)$ term is of order $O(\log |\log \alpha_{\max}|/|\log \alpha_{\max}|)$. It's noteworthy that this result holds true not only in the asymptotically 
symmetric case where $\log \alpha_0 \sim \log \alpha_1$ and $a_0 \sim a_1$ as $\alpha_{\max} \to 0$, but also in the asymmetric case where $a_0$ and~$a_1$ diverge 
with different rates, as long as $a_1e^{-a_0} \to 0$.

Note that the Schwarz--Lorden asymptotic theory operates under the assumption  of a fixed indifference zone that does not permit local alternatives, meaning that $\theta_1$ cannot approach~$\theta_0$ 
as~$c \to 0$.  In simpler terms, this theory is confined to scenarios where the width of the indifference zone $\theta_1-\theta_0$ is considerably larger than~$c^{1/2}$. 

We conclude this section by mentioning Lorden's~\citeyearpar{Lorden-AS73} paper on the properties of the one-sided (open-ended) GSLRTs
for the one-parameter exponential family. These tests reject a null hypothesis $\theta=\theta_0$ in favor of $\theta>\theta_0$ within the class of stopping times satisfying 
$\Pb_{\theta_0}(T<\infty) \leq \alpha$ for a prescribed $0<\alpha< 1/3$.

\subsection{Optimal Multistage Testing}\label{sec:multistage}

What is the fewest number of stages for which a  multistage hypothesis test can be asymptotically equivalent to an optimal fully sequential test?  \citet{Lorden83} took up this question and reached the definitive answer of needing 3 stages in general, except in a special symmetric situation  in which 2 stages are possible, described in the next section.  Here, ``needing 3 stages'' means allowing the \emph{possibility} of 3 stages, although Lorden's optimal procedures can (and do, with probability approaching~$1$) terminate earlier; see Section~\ref{sec:mult.comp}. \citet{Lorden83} shows this first in the simple vs.\ simple testing setup, and then for testing separated composite hypotheses in an exponential family. In this area again, Lorden's work was groundbreaking and formed the foundation for later, more general theoretical investigations in optimal multistage testing \citep[e.g.][]{Bartroff06b, Bartroff06, Bartroff07c, Xing23}
and in applications to clinical trial designs where the problem is sometimes known as ``sample size adjustment'' or ``re-estimation'' \citep[e.g.,][]{Bartroff08c, Bartroff08, Bartroff13}. In this literature especially, multistage procedures are often referred to as group sequential. Throughout this section, i.i.d.\ observations are assumed.

\subsubsection{Simple vs.\ Simple Testing: Multistage Competitors of the SPRT}\label{sec:mult.simp}
Beginning with the simple vs.\ simple testing setup of Section~\ref{ssec:PF} and adopting the notation there, some of Lorden's main ideas can be seen by first considering the symmetric case where the error probabilities~$\alpha_0,\alpha_1\to 0$ in such a way that
\begin{equation}\label{ENs.SPRT.symm}
\frac{\log\alpha_1^{-1}}{I_0}\sim \frac{\log\alpha_0^{-1}}{I_1}.
\end{equation} Letting $\lambda_n$ be the log-likelihood ratio statistic in \eqref{llr.sprt} and $t\to\infty$ an argument parameterizing $\alpha_0,\alpha_1\to 0$, \citet{Lorden83} begins by arguing that there is a sample size $n=n(t)\ge t$ such that $n=t+o(t)$,
\begin{equation}\label{P.bound.mult.simp}
\Pb_0(-\lambda_n<t I_0)\to 0,\quad\mbox{and}\quad \Pb_1(\lambda_n<t I_1)\to 0.
\end{equation} 
More explicitly, this is achievable by taking $n=t+\delta_t$ with
\begin{equation}\label{del.t.sqrt}
  \sqrt{t} \ll \delta_t\ll t  
\end{equation} since, assuming finite second moments $\Eb_i[\lambda_1^2]<\infty$, Chebyshev's inequality gives 
\begin{equation*}
\Pb_1(\lambda_n<t I_1)\le \frac{\var_1(\lambda_1)/n}{I_1^2(1-t/n)^2}
\end{equation*} which, ignoring constants, under \eqref{del.t.sqrt} is
\begin{equation*}
\frac{1/n}{(1-t/n)^2} = \frac{n}{\delta_t^2}= \frac{t+\delta_t}{\delta_t^2}=o(1)\quad\mbox{as $t\to\infty$.}
\end{equation*} A similar argument shows that the other probability in \eqref{P.bound.mult.simp} approaches $0$ as well.

In this symmetric situation, an optimal 2-stage competitor to the SPRT can be described in terms of $n(t)$, which is the size of the first stage with $t$ taken to be the larger of the two sides of \eqref{ENs.SPRT.symm}. Note that, for either $i=0$ or $1$, we have $t I_i \sim \log\alpha_{1-i}^{-1}$ so that $t$ is asymptotically the same as the expected stopping time of the SPRT (under either hypothesis) and the first stage $n(t)$ is of the same order but slightly larger. The procedure stops after the first stage if 
\begin{equation}\label{mult.stg.stop}
   \lambda_{n(t)} \not\in (\log\alpha_1, \log\alpha_0^{-1}), 
\end{equation} making the appropriate terminal decision. Using \eqref{P.bound.mult.simp}, the probability under the null of terminating and making the correct terminal decision after this first stage is 
\begin{equation}\label{P0.lambd.cross}
\Pb_0(\lambda_{n(t)}\le \log\alpha_1)= \Pb_0(-\lambda_{n(t)}\ge \log\alpha_1^{-1})\ge \Pb_0(-\lambda_{n(t)}\ge tI_0)\to 1,
\end{equation} with a similar argument showing that
\begin{equation}\label{P1.lambd.cross}
\Pb_1(\lambda_{n(t)}\ge \log\alpha_0^{-1}) \to 1.
\end{equation}
Otherwise, the test continues to a total sample size~$n_2$ which is that of the fixed-sample size test with error probabilities~$\alpha_0, \alpha_1$ and uses that terminal decision rule. This can be accomplished in at most $n_2\le Ct\le Cn(t)$ total observations for some constant~$C$. Thus, under the null, the total expected sample size is at most
\begin{equation*}
n(t)+Cn(t)\Pb_0(\lambda_{n(t)}> \log\alpha_1)=n(t)[1+o(1)]\sim t\sim \frac{\log\alpha_1^{-1}}{I_0},
\end{equation*} and is of the same order~$I_1^{-1}\log\alpha_0^{-1}$ under the alternative by a similar argument. By definition of the 2 stages, the procedure has type~I error probability at most $2\alpha_0$, and type~II error probability at most $2\alpha_1$, so repeating the construction with $\alpha_i/2$ replacing $\alpha_i$ ($i=0,1$) controls the error probabilities at the nominal levels and does not affect the asymptotic estimates above. Thus, this 2-stage procedure is asymptotically as efficient as the SPRT in this symmetric case.

If the asymptotic equivalence~\eqref{ENs.SPRT.symm} does not hold but we assume that 
\begin{equation}\label{rat.alpha.bdd}
\frac{\log\alpha_0^{-1}}{\log\alpha_1^{-1}}\quad\mbox{is bounded away from $0$ and $\infty$,}  
\end{equation} Lorden shows that no 2-stage test can be asymptotically optimal, itself a nontrivial result that we discuss in the next section. For this case Lorden gives a 3-stage procedure that is a slight modification of the one above.  Letting $t_1$ and $t_2$ be the left- and right-hand sides of \eqref{ENs.SPRT.symm}, respectively, the first stage of the procedure is of size $\min\{n(t_1), n(t_2)\}$, and the second stage (if needed) brings the total sample size to $\max\{n(t_1), n(t_2)\}$, both using the stopping rule~\eqref{mult.stg.stop} and corresponding decision rule.  If not stopped by the second stage, a third stage brings the total sample size to that of the fixed-sample size with error probabilities $\alpha_i$, which is $\le C \max\{n(t_1), n(t_2)\}$ as above, and uses that terminal decision rule.  Since $n(t_{i+1})\sim t_{i+1}\sim \log \alpha_{1-i}^{-1}/I_i$ for both $i=0,1$, and \eqref{P0.lambd.cross} and \eqref{P1.lambd.cross} hold for $n(t_1)$ and $n(t_2)$, respectively, the expected sample size of this 3-stage procedure is asymptotically equal to the corresponding side of \eqref{ENs.SPRT.symm}, and is thus minimized under both the null and alternative.

\subsubsection{The Necessity of 3 Stages}\label{sec:3.need}

Continuing with the simple vs.\ simple testing setup of the previous section, Lorden's~\citeyearpar[][Corollary~1]{Lorden83} result  mentioned above that, in the absence of symmetry~\eqref{ENs.SPRT.symm}, 3 stages are necessary for asymptotic optimality, is far from obvious since it may seem that the first 2 stages of the 3 stage procedure defined above would suffice.  That is, why is it that a first stage of $\min\{n(t_1), n(t_2)\}$ and (if needed) a second stage giving total sample size $\max\{n(t_1), n(t_2)\}$ would not be optimal? One clue may  be that, if that were true, then the same reasoning would seem to imply that a single-stage test could be optimal under symmetry~\eqref{ENs.SPRT.symm}, which is known to not hold.  More generally, Lorden provides the following general result about asymptotically optimal $k$-stage ($k\ge 2$) tests: that their expected sample size after $k-1$ stages must be asymptotically the same as after $k$ stages. In other words, the final stage of an asymptotically optimal multistage test is asymptotically negligible in size, but \emph{necessary}. In what follows let $I(f,g)$ denote the information number for arbitrary densities $f,g$.

\begin{theorem}[\citet{Lorden83}, Theorem~3]\label{thm:lord.k.stg}
For testing $f_0$ vs.\ $f_1$ in the setup of Section~\ref{ssec:PF}, let $N$ denote the sample size of a $k$-stage ($k\ge 2$) test with error probabilities $\alpha_0$ and $\alpha_1$, and let $M$ be the total sample size of this test after $k-1$ stages. If $N$ is asymptotically optimal as $\alpha_0,\alpha_1\to 0$ and $g$ is a density distinct from $f_0$ such that 
\begin{equation*}
\frac{\log\alpha_1^{-1}}{\log\alpha_0^{-1}} \ge Q> \frac{I(g,f_1)}{I(g,f_0)}
\end{equation*} for some $Q>0$ as $\alpha_0,\alpha_1\to 0$, then
$$M\to \frac{\log\alpha_0^{-1}}{I(g,f_0)} \quad\mbox{in $g$-probability, and}\quad  \Eb_g [M] \sim \frac{\log\alpha_0^{-1}}{I(g,f_0)}\sim \Eb_g[N]\quad\mbox{as $\alpha_0\to 0$.}$$ 
\end{theorem}

Lorden's proof of this theorem is technical and requires detailed upper bounds on the \emph{conditional} error probabilities after the $(k-1)$st stage; that is, the probabilities of test error given the first $M$ observations. Roughly speaking, showing that these error probabilities are small shows that their corresponding sample size~$M$ must be large, so large in fact that it is asymptotically equivalent to its maximum value~$N$.

\citet[][Corollary~1]{Lorden83} then uses Theorem~\ref{thm:lord.k.stg} to show that there is an asymptotically optimal 2-stage test \emph{if and only if} the symmetry condition~\eqref{ENs.SPRT.symm} holds, with the construction of the 2-stage test above providing the ``if'' argument.  For the converse, applying Theorem~\ref{thm:lord.k.stg} with $g=f_1$ shows that the first stage of an optimal 2-stage test must be asymptotic to $(\log\alpha_0^{-1})/I_1$. After reversing the roles of $f_0$ and $f_1$ in the theorem and applying it again with $g=f_0$, it also shows that the first stage must be asymptotic to $(\log\alpha_1^{-1})/I_0$, establishing symmetry~\eqref{ENs.SPRT.symm}.

\subsubsection{Composite Hypotheses}\label{sec:mult.comp}
For testing separated hypotheses $\theta\le\theta_0$ vs.\ $\theta\ge \theta_1>\theta_0$ about the 1-dimensional parameter~$\theta$ of an exponential family, \citet[][Section~3]{Lorden83} constructs an asymptotically optimal 3-stage test utilizing a description of the optimal stopping boundary related to Schwarz's~\citeyearpar{Schwarz-AMS1962} study of Bayes asymptotic shapes for fully sequential tests, described in Section~\ref{ssec:GLR}.  Let $n(\theta)$ denote the expected sample size to Schwarz's boundary under $\theta$. Lorden's test utilizes the ``worst case'' competing parameter value~$\theta^*\in(\theta_0,\theta_1)$ which maximizes the expected sample size~$n(\theta^*)=\max_\theta n(\theta)\equiv n^\star$. The first stage size of Lorden's procedure is a fixed fraction of $n^\star$. If the procedure does not stop after the first stage, utilizing Schwarz's boundary, the second stage brings the total sample size to
$\min\{n^\star,(1+\varepsilon)n(\widehat{\theta})\}$, where $\widehat{\theta}$ is the MLE of $\theta$ from the first stage data and $\varepsilon\searrow 0$ is a chosen sequence. Finally, if needed, the third stage brings the total sample size up to $n^\star$. Under \eqref{rat.alpha.bdd}, \citet[][Theorem~1]{Lorden83} proves that this test asymptotically minimizes the expected sample size to first order, not just for $\theta$ in the hypotheses but uniformly in $\theta$ over any interval in the parameter space containing $[\theta_0,\theta_1]$. The first order term is of order $\log\alpha_i^{-1}$, as above, and the second order term is of order $O(((\log\alpha_i^{-1})\log\log\alpha_i^{-1})^{1/2})$, $i=0,1$. 

These results were extended to asymptotically optimal 3-stage tests of multidimensional parameters in \citet{Bartroff06b} and \citet{Bartroff08c}, and more general multidimensional composite hypotheses in \citet{Bartroff08}. On the other hand, Lorden's procedures were generalized to optimal $k$-stage tests, for arbitrary $k\ge 3$, in \citet{Bartroff06,Bartroff07c}.

Regarding the necessity of 3 stages in this composite hypothesis setting,  \citet[][Corollary~2]{Lorden83} proves that, under \eqref{rat.alpha.bdd}, 3 stages are necessary (and sufficient, by his own procedure) for asymptotic optimality at more than 3 values of $\theta$, and so certainly for asymptotic optimality over an interval of $\theta$ values, as in Lorden's result.  An interesting detail that shows this result to be best possible is that an optimal 2-stage test \emph{can} be constructed at 3 values of $\theta$ if the  special symmetry condition $I(\theta',\theta_0)I(\theta_0,\theta_1) = I(\theta',\theta_1)I(\theta_1,\theta_0)$ holds for some $\theta'\ne \theta_0,\theta_1$. Then a 2-stage procedure similar to the one described in Section~\ref{sec:mult.simp} that uses second stage total sample size of $\log\alpha_0^{-1}/I(\theta',\theta_0)$ will be optimal at the 3 values $\theta=\theta'$, $\theta_0$, and $\theta_1$.

\section{Sequential Changepoint Detection: Lorden's Minimax Change Detection Theory} \label{sec:CPD}

In numerous practical applications, the observed process undergoes an abrupt change in statistical properties at an unknown point in time. Examples encompass aerospace navigation and flight systems integrity monitoring, cyber-security, identification of terrorist activity, industrial monitoring, air pollution monitoring, radar, sonar, and electrooptics surveillance systems. Consequently, this problem has garnered interest from many practitioners for some time.

In classical quickest changepoint detection, the objective is to detect changes in the distribution as swiftly as possible, thereby minimizing the expected delay to detection assuming the change is in effect.

More specifically, the changepoint problem posits that one obtains a series of observations $X_1,X_2,\dots$ such that, for some value $\nu$, $\nu \in \Zbb_+ =\{0,1, 2, \dots\}$  
(the changepoint), $X_1,X_2,\dots,X_{\nu}$ have one distribution and $X_{\nu+1}, X_{\nu+2},\dots$ have another distribution. The changepoint~$\nu$ is unknown, and the sequence 
$\{X_n\}_{n\ge 1}$ is being monitored for detecting a change. A sequential detection procedure is a stopping time~$T$ with respect to the~$X$s, 
so that after observing $X_1,X_2,\dots,X_T$ it is declared that a change is in effect. That is, $T$ is an integer-valued random variable, such that the event $\{T = n\}$
belongs to the sigma-algebra  $\Fc_{n}=\sigma(X_1,\dots,X_n)$ generated by observations $X_1,\dots,X_n$.

Historically, the field of changepoint detection began to take shape in the 1920s to 1930s, spurred by considerations in quality control. Shewhart's charts were particularly influential during this period~\citep{shewhart-book31}. However, optimal and nearly optimal sequential detection procedures didn't come into prominence until much later, in the 1950s to 1970s, following the advent of Sequential Analysis~\citep{wald47}. The concepts initiated by Shewhart and Wald laid the foundation for extensive research into sequential changepoint detection.

The desire to detect the change quickly often leads to being ``trigger-happy,'' which, on one hand, results in an unacceptably high false alarm rate -- terminating the process prematurely before a real change has occurred. On the other hand, attempting to avoid false alarms too strenuously causes a long delay between the true change point and its detection. Thus, the essence of the problem lies in achieving a tradeoff between two conflicting performance measures -- the loss associated with the delay in detecting a true change and that associated with raising a false alarm. An efficient detection procedure is expected to minimize the average loss associated with the detection delay, while subject to a constraint on the loss associated with false alarms, or vice versa.

Let $p_\nu(\Xb^n)= p(X_1,\dots,X_n | \nu)$ denote the joint probability density of the sample $\Xb^n=(X_1,\dots, X_n)$ when the changepoint~$\nu$ is fixed ($0 \leq \nu <\infty$) and 
$p_\infty(\Xb^n)$ the joint density when $\nu=\infty$, i.e., when there is never a change. Let $\Pb_\nu, \Pb_\infty$ and $\Eb_\nu, \Eb_\infty$ denote the corresponding 
probability measures and expectations. 
Assume  that the observations $\{X_n\}_{n\ge1}$ are independent and such that $X_1,\ldots,X_{\nu}$ are each distributed according to a common (pre-change) density~$f_0(x)$, 
while $X_{\nu+1},X_{\nu+2},\ldots$ each follows a common (post-change) density $f_1(x)$.  Hence, the  model can be represented as
\begin{equation}\label{iidmodel}
p_{\nu}(\Xb^n) = 
\begin{cases}
\prod_{t=1}^{\nu} f_0(X_t) \times \prod_{t=\nu+1}^n f_1(X_) & \text{for}~n \geq \nu+1 
\\
p_\infty(\Xb^n) =\prod_{t=1}^{n} f_0(X_t) & \text{for}~ 1 \leq n \leq \nu 
\end{cases} .
\end{equation}
Note that we assume that $X_{\nu}$ is the last pre-change observation, which is different from many publications (including Lorden's) where it is assumed that 
$X_\nu$ is the first post-change observation. The diagram below illustrates this case
\[
    \underbrace{X_1, \cdots , X_\nu }_{\text{i.i.d., $f_0$}}, ~
    \underbrace{X_{\nu + 1}, X_{\nu + 2}, \cdots}_{\text{i.i.d., $f_1$}}.
\]

Denote by $\Hyp_\infty: \nu=\infty$ the hypothesis that the change never occurs and by $\Hyp_{\nu}$  the hypothesis that the change occurs at time $0 \leq \nu<\infty$. 
Let $Z_t=  \log [f_1(X_t)/ f_0(X_t)]$ denote the LLR for the $t$-th observation $X_t$.

We now introduce the CUMULATIVE SUM (CUSUM) detection procedure, which was first proposed by \citet{page-bka54}.  The changepoint detection problem can be
viewed as a problem of testing two hypotheses: $\Hyp_\nu$ that the change occurs at 
a fixed point $0 \leq \nu <\infty$ against the alternative $\Hyp_\infty:\nu=\infty$ that the change never occurs. The LLR between these hypotheses
 is $\lambda_n^\nu = \sum_{t=\nu+1}^n Z_t$ for $\nu<n$ and~$0$ for $\nu\geq n$. Since the hypothesis~$\Hyp_\nu$ is composite, we may employ the GLR 
 approach, maximizing the LLR $\lambda_n^\nu$ over~$\nu$, to obtain the log-GLR statistic:
\begin{equation}\label{CUSUMstat}
W_n =  \max_{\nu \geq 0}\sum_{t=\nu+1}^{n}Z_t,
\end{equation}
which follows the recursion
\begin{equation}\label{CUSUMstat1}
W_n = \brc{W_{n-1} + Z_n}^+, \quad  n \ge 1,~~ W_0=0.
\end{equation}
This statistic is called the CUSUM statistic. Page's CUSUM procedure is the first time $n \geq 1$ such that the CUSUM statistic $W_n$ exceeds a positive threshold $a$:
\begin{equation}\label{CUSUM_ST}
T_a = \inf\{n \ge 1 : W_n \ge a\}.
\end{equation}

\citet{page-bka54} proposed measuring the risk due to a false alarm by the mean time to false alarm $\Eb_\infty [T]$
 and the risk associated with a true change detection by the mean time to detection $\Eb_0 [T]$ when the change occurs at the very beginning. 
 These are commonly known as the Average Run Length (ARL). Page also analyzed the CUSUM procedure defined by equations
 \eqref{CUSUMstat}--\eqref{CUSUM_ST} using these operating characteristics.

While the false alarm rate is reasonable to measure by the ARL to false alarm $\ARLFA(T)=\Eb_\infty[T]$, the risk due to a true change detection is better measured by the conditional expected delay to detection $\Eb_\nu[T-\nu | T > \nu]$ for any possible change point $\nu \in \Zbb_+$, rather than by the ARL to detection $\Eb_0[T]$.
Ideally, a good detection procedure should guarantee small values of the expected detection delay for all change points  $\nu\in \Zbb_+$ when $\ARLFA(T)$
 is set at a certain level. However, if the false alarm risk is measured in terms of the ARL to false alarm, i.e., it is required that $\ARLFA(T)\geq \gamma$ for some $\gamma \geq 1$,
 then a procedure that minimizes the conditional expected delay to detection $\Eb_\nu[T-\nu | T > \nu]$ uniformly over all~$\nu$
 does not exist. For this reason, we must resort to different optimality criteria, such as Bayesian and minimax criteria.

The minimax approach posits that the changepoint is an unknown not necessarily random number. Even if it is random its distribution is unknown. 

\citet{lorden-ams71} was the first who addressed the minimax change detection problem and developed the first minimax theory. He proposed to measure the false alarm risk by the 
ARL to false alarm $\ARLFA(T)=\Eb_{\infty}[T]$, i.e., to consider the class of change detection procedures  
$ \class(\gamma)= \set{T: \ARLFA(T) \geq \gamma}$ for some $\gamma \geq 1$,
 and the risk associated with detection delay by the worst-case expected detection delay 
\begin{equation}\label{eq:SADD-Lorden-def}
\ESADD(T) = \sup_{0 \leq \nu<\infty}\biggl\{\esssup \Eb_{\nu}[(T-\nu)^+| X_1,\dots, X_\nu]\biggr\}.
\end{equation}
In other words, the conditional expected detection delay is maximized over all possible trajectories $(X_1,\dots,X_\nu)$ up to the changepoint 
and then over the changepoint $\nu$. 

Lorden's minimax criterion is
\begin{equation*}
    \inf_{T} \sup_{\nu \geq 0} \esssup_\omega \Eb_{\nu}[T - \nu \mid T > \nu, \Fc_\nu] \quad \text{subject to}~ \ARLFA(T) \geq \gamma,
\end{equation*}
i.e.,  Lorden's minimax optimization problem seeks to
\begin{equation}\label{eq:Lorden-minimax-problem}
\text{Find $T_{\mathrm{opt}}\in\class(\gamma)$ such that $\ESADD(T_{\mathrm{opt}})= \inf_{T\in\class(\gamma)}\ESADD(T)$ for every $\gamma \geq 1$} .
\end{equation}

\citet{lorden-ams71} demonstrated that Page's CUSUM procedure achieves first-order asymptotic minimax optimality as $\gamma$ approaches infinity. This groundbreaking finding marked the initial optimality result in the minimax change detection problem. Given the significance of this outcome and the widespread adoption of Lorden's minimax criterion not only within statistical circles but also across various practical domains, we proceed to provide further elaboration.

To establish the asymptotic optimality of Page's CUSUM procedure, Lorden employs an intriguing method that permits the utilization of one-sided hypothesis tests to assess a collection of change detection procedures, among them Page's method.  Let $\tau=\tau(\alpha)$ be a stopping time with respect to $X_1,X_2,\ldots$ such that 
\begin{equation}\label{Pinftyalpha}
\Pb_\infty(\tau<\infty) \le \alpha,
\end{equation}
where $\alpha\in(0,1)$. For $k=0,1,2,\dots$ define the stopping time $\tau_k$ obtained by applying $\tau$ to the sequence $X_{k+1}, X_{k+2}, \dots$ and let 
$\tau^*=\min_{k\ge 0}(\tau_k+k)$.

The subsequent theorem, resembling Theorem~2 in \citet{lorden-ams71}, empowers the construction of nearly optimal change detection procedures and facilitates the demonstration of the near optimality of the CUSUM procedure. It's important to recall that $\Pb_\infty$ denotes the distribution characterized by the density $f_0(x)$, while $\Pb_0$ corresponds to the distribution with density $f_1(x)$.

\begin{theorem}\label{th:L2}
The random variable $\tau^*$ is a stopping time with respect to $X_1,X_2, \dots$ and if condition \eqref{Pinftyalpha} is satisfied, then the following two inequalities hold:
\begin{equation}\label{ARL}
\Eb_\infty[\tau^*] \ge  1/\alpha
\end{equation}
and
\begin{equation}\label{ESS}
\Eb_0[\tau^*] \le  \Eb_0[\tau].
\end{equation}
\end{theorem}

The cumulative LLR for the sample $(X_{k+1},\dots,X_n)$ is 
$
\lambda_n^k = \sum_{t=k+1}^n Z_t. 
$
Let $\tau(\alpha)=\inf\set{n \ge 1: \lambda_n^0 \ge |\log\alpha|}$ denote the stopping time of the one-sided SPRT for testing $f_0$ {\it versus} $f_1$ with threshold $|\log\alpha|$. Then $\Pb_\infty(\tau(\alpha)<\infty) \le \alpha$, so condition \eqref{Pinftyalpha} holds. If the Kullback-Leibler information number
$I =\Eb_0[Z_{1}]$ is positive and finite, then it is well-known that
\[
\Eb_0[\tau(\alpha)] = \frac{|\log \alpha|}{I}(1+o(1)) \quad \text{as}~ \alpha\to 0.
\]

Next, note that the CUSUM statistic defined in \eqref{CUSUMstat} is the maximum of $\lambda_n^k$ over $k \ge 0$, so the stopping time of the CUSUM procedure 
\eqref{CUSUMstat1} can obviously be written as $T_a=\min_{k\ge 0} \{\tau_k(\alpha)+k\}\equiv \tau^*$ for $a= a_\alpha=|\log\alpha|$, where
\[
\tau_k(\alpha) = \inf\set{n \ge 1: \lambda_{k+n}^k \ge |\log\alpha|}.
\]
It follows from Theorem~\ref{th:L2} that setting $\alpha=\gamma^{-1}$ gives $\Eb_\infty[T_{a_\gamma} ] \ge \gamma$, so $T_{a_\gamma} \in \class(\gamma)$, and
\[
\ESADD(T_{a_\gamma}) \equiv \Eb_0[T_{a_\gamma}] = \frac{\log \gamma}{I}(1+o(1)) \quad \text{as}~ \gamma\to \infty.
\]

To complete the proof of the first-order asymptotic optimality of the CUSUM procedure with threshold $a=a_\gamma= \log \gamma$ it suffices to establish that this is the best one can do, 
i.e., to prove the asymptotic lower bound
\begin{equation}\label{LB}
\inf_{T\in\class(\gamma)} \ESADD(T) \ge \frac{\log \gamma}{I}(1+o(1)) \quad \text{as}~ \gamma\to \infty,
\end{equation}
which also yields 
\[
\inf_{T\in\class(\gamma)} \ESADD(T) \sim \frac{\log \gamma}{I} \sim \ESADD(T_{a_\gamma}) \quad \text{as}~ \gamma\to \infty.
\]
Theorem 3  of \citet{lorden-ams71} establishes this fact using a rather sophisticated argument. Note, however, that \citet{LaiIEEE98} established the lower bound \eqref{LB} in a general non-i.i.d. case, assuming that $n^{-1}\lambda_{\nu+n}^\nu$ converges to a positive and finite number $I$ as $n\to\infty$, under a certain additional condition.
In the i.i.d.\ case, by the SLLN $n^{-1}\lambda_{\nu+n}^\nu$ converges to the Kullback--Leibler information number $I$ almost surely under $\Pb_\nu$. This implies that as 
$M\to \infty$ for all $\varepsilon >0$
\begin{equation}\label{Pmax0}
\sup_{\nu \ge 0} \Pb_\nu\set{\frac{1}{M} \max_{0\le n \le M}\lambda_{\nu+n}^\nu \ge (1+\varepsilon)I} = \Pb_0\set{\frac{1}{M}\max_{0\le n \le M}\lambda_{n}^0 \ge (1+\varepsilon)I} \to 0.
\end{equation}
Using \eqref{Pmax0}, the lower bound \eqref{LB} can be obtained from Theorem~1 in \citet{LaiIEEE98}.

To handle a composite parametric post-change hypothesis, which is typical in many applications,  let $f_{\theta}(x)$ be the post-change density, where $\theta\in \Theta$. 
Denote $Z_{n}(\theta) =\log[f_{\theta}(X_n)/f_0(X_n)]$. Then, inequality \eqref{ESS} in Theorem~\ref{th:L2} holds for expectation $\Eb_\theta[\tau^*]$. 
Additionally, assuming that the Kullback-Leibler information number $I(\theta) =\Eb_\theta[Z_{1}(\theta)]$ is positive and finite, then asymptotic lower bound \eqref{LB} holds
with $I(\theta)$, i.e.,
\begin{equation}\label{LBtheta}
\inf_{T\in\class(\gamma)} \ESADD_\theta(T) \ge \frac{\log \gamma}{I(\theta)}(1+o(1)) \quad \text{as}~ \gamma\to \infty,
\end{equation}
where $\ESADD_\theta(T) = \sup_{0 \leq \nu<\infty} \esssup \Eb_{\nu, \theta}[(T-\nu)^+| \Fc_\nu]$ and $\Eb_{\nu,\theta}$ is the expectation under $\Pb_{\nu,\theta}$
when the change occurs at $\nu$ with the post-change density $f_\theta$.

 \citet{lorden-ams71} addressed the composite hypothesis for the exponential family \eqref{Expfam} with $f_0=f_{\theta=0}$, i.e.,
 \begin{equation*}
\frac{f_\theta(X_n) }{f_{0}(X_n)} = \exp\set{\theta X_n -b(\theta)}, \quad \theta\in \Theta, \quad n=1,2, \dots
\end{equation*}
where~$b(\theta)$ is a convex and infinitely differentiable function on the natural parameter space $\Theta$, $b(0)=0$. Let $\widetilde \Theta= \Theta-0$.

In order to find asymptotically optimal procedures by applying Theorem~\ref{th:L2} along with inequality \eqref{LBtheta} we need to determine stopping times, 
$\tau(\gamma)\in\class(\gamma)$, such that
\begin{equation}\label{Probupper}
\Pb_0(\tau(\gamma) < \infty) \leq 1/\gamma \quad \text{for} ~ \gamma >0
\end{equation}
and
\begin{equation}\label{Exptheta}
\Eb_\theta[\tau(\gamma)] = \frac{\log \gamma}{I(\theta)}(1+o(1)) \quad \text{as}~ \gamma\to \infty \quad \text{for all}~ \theta\in\widetilde\Theta,
\end{equation}
where $I(\theta) =  \theta \dotb{b}(\theta)  -b(\theta)$.

The LLR for the sample $(X_{k+1},\dots,X_n)$ is
\[
\lambda_n^k(\theta) =: \log \brcs{\prod_{t=k+1}^n\frac{f_\theta(X_t)}{f_0(X_t)}} = \theta S_n^k - (n-k) b(\theta) , 
\]
where $S_n^k=X_{k+1}+ \cdots + X_n$. Define the GLR one-sided test 
\[
\tau(h)= \inf\set{n\ge 1: \sup_{\theta\ge |\theta_1|} \brcs{ \theta S_n^0 - n b(\theta)} >h(\gamma)},
\]
where $\theta_1$ may be either a fixed value if the alternative hypothesis is $\theta \le -\theta_1$ or $\theta \ge \theta_1$ or $\theta_1(\gamma) \to 0$ as $\gamma\to \infty$ if the 
hypothesis is $\theta\neq 0$. Lorden demonstrates that 
\begin{equation} \label{ProbtauUB}
\Pb_0(\tau(h) < \infty) \le \exp\set{-h(\gamma)} \brcs{1+ \frac{h(\gamma)}{\min(I(\theta_1),I(-\theta_1)}},
\end{equation}
so $h(\gamma)$ can be selected so that $h(\gamma) \sim \log \gamma$ as $\gamma \to \infty$. Hence, \eqref{Probupper} and \eqref{Exptheta} hold.
Applying $\tau(h)$ to $X_{k+1}, X_{k+2},\dots$ we obtain the stopping time $\tau_k(h)$, 
so that $\tau^*(h)=\min_{k\ge 0}(\tau_k+k)$ is the stopping time of the GLR CUSUM procedure,
\[
\tau^*(h) = \inf\set{n \ge 1: \max_{ 0 \le \nu \le n}  \sup_{\theta\ge |\theta_1|} \brcs{ \theta S_n^\nu - (n-\nu) b(\theta)} >h(\gamma)}.
\]
Thus, the GLR CUSUM procedure is asymptotically first-order minimax. 

The inequality \eqref{ProbtauUB} is usually overly pessimistic. A much better result gives the approximation
$
\Pb_0(\tau(h) < \infty) \approx \sqrt{h(\gamma)}\exp\set{-h(\gamma)} C,
$
which follows from \eqref{PEGSLRTlorden}. However, the latter one does not guarantee the inequality $\Pb_0(\tau(h) < \infty) \le \gamma^{-1}$, and therefore, the 
inequality $\Eb_\infty[\tau^*(h)] \ge \gamma$.

Later, \citet{LordenPollak-as05,LordenPollak-SQA08} proposed adaptive Shiryaev-Roberts and CUSUM procedures that utilize one-step delayed estimators of unknown post-change parameters $\theta$.  In these procedures, an estimate $\hat{\theta}_{n-1}(X_1,\dots,X_{n-1})$ is used after observing the sample of size $n$,
similar to the Robbins--Siegmund one-sided adaptive SPRT; see \citet{RobbinsSiegmund-Berkeley70,RobbinsSiegmund-AS74}. They compared the performance of these adaptive procedures with that of the mixture-based Shiryaev-Roberts procedure. Notably, these adaptive procedures are computationally simpler than the GLR CUSUM procedure.

We conclude with some remarks on later, related developments.

\vspace{3mm}

{\bf REMARKS}
\vspace{2mm}

1.   Fifteen years later, \citet{MoustakidesAS86} advanced Lorden's asymptotic theory by demonstrating, using optimal stopping theory, that the CUSUM procedure is 
strictly optimal for any ARL to false alarm $\gamma\ge 1$ if the threshold $a=a(\gamma)$ is chosen such that $\ARLFA(T_a)= \gamma$. 

2. \citet{shiryaev-rms96} showed that the CUSUM procedure is strictly optimal in the continuous-time scenario for detecting the change in the mean of the Wiener process according to Lorden's minimax criterion.

3.  \citet{PollakAS85} introduced a distinct minimax criterion aimed at minimizing the supremum expected detection delay  $\sup_{\nu\ge 0} \Eb[T-\nu |T > \nu]$. 
Additionally, Pollak proposed a modification of the conventional Shiryaev--Roberts (SR) procedure known as the SRP procedure, which initiates from a randomly distributed 
point following the quasi-stationary distribution of the SR statistic. He proved that  this procedure is third-order asymptotically minimax, minimizing 
$\sup_{\nu\ge 0} \Eb[T-\nu |T > \nu]$ to within $o(1)$ as $\gamma \to \infty$ within the class $\class(\gamma)$.

4. \citet{tartakovskypolpolunch-tpa11} proved that the specially designed SR-$r$ procedure that starts from a fixed point $r=r(\gamma)$ is third-order asymptotically optimal 
with respect to Pollak's measure $\sup_{\nu\ge 0} \Eb[T-\nu |T > \nu]$ within the class $\class(\gamma)$ as $\gamma \to \infty$.

5. \citet{PolunTartakovskyAS10} demonstrated that the specially designed SR-$r$ procedure, which commences from a predetermined point $r=r(\gamma)$, is strictly optimal with respect to Pollak's measure $\sup_{\nu\ge 0} \Eb[T-\nu | T > \nu]$ within the class $\class(\gamma)$ for a specific model.

6. \citet{PollakTartakovsky-SS09} proved strict optimality of the repeated SR procedure that starts from zero in the problem of detecting distant changes.

7. \citet{MoustPolTarCS09} conducted a thorough comparison of CUSUM and SR procedures, demonstrating that CUSUM outperforms SR in terms of the conditional expected detection delay $\Eb_\nu[T-\nu | T> \nu]$ for relatively small values of the change point $\nu$. However, SR proves to be more effective than CUSUM for relatively large $\nu$.

\section*{Acknowledgements}
We express our gratitude to the reviewers,  the associate editor, and the editor for their valuable comments, which greatly enhanced the quality of the paper.  We also thank Caltech for the use of the photos in Figure~\ref{fig:lorden_pics}.

AT: I am grateful to Gary Lorden for multiple helpful and insightful conversations starting in 1993 and Gary's many papers that we have discussed in this article. Gary's work meaningfully influenced my research, from 1977 on.

JB: I am lucky to be able to call Gary Lorden not only my PhD advisor, but also a mentor and friend.  Gary had a profound effect on my life, both within and  ``beyond the boundaries'' of academics.


\end{document}